\documentclass[11pt]{article}
\usepackage{indentfirst}
\usepackage{amsfonts}
\usepackage{amsmath}
\usepackage{epsfig}
\usepackage{amssymb} % misc. symbols
\usepackage{enumerate} % misc. symbols
\usepackage{mathrsfs}
\usepackage{graphicx}
\usepackage{hhline}
\usepackage{array}
\usepackage{cite}
\usepackage{wrapfig}%does not work in enumerate environment
\usepackage{float}%put figure at a reqired place
\usepackage{color}

\def\qed{\hspace*{\fill}$\Box$\medskip}
\newcommand{\proof}{\noindent{\bf Proof.}\ \ }
\newcommand{\round}[1]{\mbox{\textcircled{\scriptsize{$#1$}}}}

\def\rsq{\hspace*{\fill}$\square$\medskip}
\newcommand{\bb}{$\mathcal{B}$\ }

\newtheorem{theorem}{Theorem}
\newtheorem{lemma}[theorem]{Lemma}
\newtheorem{defi}[theorem]{Definition}

\newenvironment{definition}{\begin{defi}\rm }{\end{defi}}
\newtheorem{ex}{Example}[section]
\newtheorem{rem}{Remark}[section]

\begin{document}

\begin{center}
{\mathversion{bold}\Large\bf A characterization of $L(2,1)$-labeling number for trees with maximum degree~3}\footnote{Research supported partially by NSFC (Nos. 11401535 and 11371328); Faculty Research Grant of Hong Kong Baptist University.}
\end{center}

\begin{center}Dong Chen$^{a}$,  Wai Chee Shiu$^{b}$\footnote{The
corresponding author. Email: wcshiu@hkbu.edu.hk}, Qiaojun Shu$^{c}$, Pak Kiu Sun$^{b}$, Weifan Wang$^{d}$

\vspace{6pt}

\small{$^{a}$School of Mathematics and Statistics, Lanzhou
University, Lanzhou,}\\ \small{Gansu 730000, China}\\
\small{$^b$Department of Mathematics, Hong Kong Baptist University,}\\
\small{Kowloon Tong, Hong Kong,
China} \\
\small{$^d$Department of Mathematics, Zhejiang Normal University,\\ Jinhua, 321004, China.}

\end{center}

%
%\author[Shu,Chen]{Dong Chen},\ead{258499@qq.com}
%\author[Shiu]{Wai Chee Shiu\thanksref{correspond}}, \ead{wcshiu@hkbu.edu.hk}
%\author[Shu]{Qiaojun Shu}, \ead{shuqiaojun@163.com}
%\author[Shiu]{Pak Kiu Sun}, \ead{lionel@hkbu.edu.hk}
%\author[Wang]{Weifan Wang} \ead{wwf@zjnu.cn}
%\address[Shu]{Department of Mathematical Sciences, Soochow University,\\ Suzhou 215006, China.}
%\address[Shiu]{Department of Mathematics,
%Hong Kong Baptist University,\\
%224 Waterloo Road, Kowloon Tong, Hong Kong, China.}
%\thanks[correspond]{Corresponding author}
%\address[Chen]{Xingzhi College, Zhejiang Normal University,\\ Jinhua, 321004, China.}
%\address[Wang]{Department of Mathematics, Zhejiang Normal University,\\ Jinhua, 321004, China.}
\begin{abstract}
An $L(2, 1)$-labeling of a graph is an assignment of nonnegative integers to the vertices of $G$ such that adjacent vertices receive numbers differed by at least 2, and vertices at distance 2 are assigned distinct numbers. The $L(2,1)$-labeling number is the minimum range of labels over all such labeling. It was shown by Griggs and Yeh [Labelling graphs with a condition at distance 2, {\it SIAM J. Discrete Math.} 5(1992), 586-595] that the $L(2,1)$-labeling number of a tree is either $\Delta+1$ or $\Delta+2$. In this paper, we give a complete characterization of $L(2,1)$-labeling number for trees with maximum degree 3.

\def\MSC{\par\leavevmode\hbox{\em AMS 2010 MSC:\ }}%
\MSC 05C15.
\end{abstract}

%%%%%%%%%%%%%%%%%%%%%%%%%%%%%%%%%%%%%%%%%%%%%%%%%%%%%%%%%%%%%%%%

\section{Introduction}

An $L(2, 1)$-labeling $f$ of a graph $G$ is a function from the vertex set $V(G)$ to the set of nonnegative integers such that $|f(x)-f(y)|\ge 2$ if $x$ and $y$ are adjacent and $|f(x)-f(y)|\ge 1$ if $x$ and $y$ are at distance 2, where $x,y\in V(G)$. A $k$-$L(2,1)$-labeling of a graph is an $L(2,1)$-labeling with image at most $k$. The $L(2,1)$-labeling number of $G$, denoted by $\lambda(G)$,  is the smallest $k$ such that $G$ has a $k$-$L(2,1)$-labeling.

The $L(2,1)$-labeling of a graph arose from a variation of the Frequency Channel Assignment problem introduced by Hale \cite{HALE1980}. This subject has been studied rather extensively \cite{CHANG1996,CHANG2000,GEORGES2000,GEORGES1994,GONCALVES2005,
GRIGGS1992,HEUVEL2003,KRAL2003,SAKAI1994,WHITTLESEY1995}. It is obvious that $\lambda(G)\ge \Delta+1$ for any graph $G$, where $\Delta$ is the maximum degree of  $G$. For the upper bound of $\lambda(G)$, Griggs and Yeh \cite{GRIGGS1992}  proved $\lambda(G)\le \Delta^2+2\Delta$. Moreover, they conjectured that $\lambda(G)\le \Delta^2$ for any graph $G$ with $\Delta\ge2$ and they confirmed the conjecture for a few classes of graphs such as paths, cycles, trees, graphs with diameter 2, etc.  In 1996, Chang and Kuo \cite{CHANG1996} proved that $\lambda(G)\le \Delta^2+\Delta$ for any graph $G$. Kr\'{a}l and \u{S}krekovski \cite{KRAL2003} improved this bound by showing that $\lambda(G)\le \Delta^2 + \Delta-1$. In 2005, Daniel Gon\c{c}alves \cite{GONCALVES2005} proved that $\lambda(G)\le \Delta^2 + \Delta-2$ and this remains the best upper bound until now.

Let $T$ be a tree, Griggs and Yeh \cite{GRIGGS1992} proved that $\lambda(T)$ is either $\Delta+1$ or $\Delta+2$. It is easy to see that there exist infinitely many trees $T$ such that $\lambda(T)=\Delta+1$ or $\lambda(T)=\Delta+2$. Wang \cite{WWF2006} and Zhai \cite{ZMQ2012} established some sufficient conditions for $\lambda(T)=\Delta+1$. Naturally, it is very interesting to give a characterization for trees to have different $L(2,1)$-labeling numbers. In this paper, we give a solution for the case of $\Delta=3$.

\section{Structural analysis}

For a tree $T$, a vertex of degree $k$ is called a $k$-vertex. Let $V_k(T)$ be the set of $k$-vertices in $T$. A vertex is called {\it major} if $d(v)=\Delta$; {\it minor} if $d(v)< \Delta$; a {\it leaf} if $d(v)=1$, and a {\it handle} if $d(v)> 1$ and $v$ is adjacent to at least $d(v)-1$ leaves. A {\it $k$-handle} is a handle of degree $k$. A path $P=ux_1x_2\cdots x_kv$ is called a {\it $k$-$uv$-chain} ({\it $k$-chain} or {\it $uv$-chain}) if $d(u)\neq 2$, $d(v)\neq 2$ and $d(x_i)=2$ for all $i=1,2,\dots, k$. In particular, $P=uv$ is defined as a 0-chain. We denote a directed $k$-$uv$-chain by $[u(k)v]$ (or $[uv]$ for short) if it is a directed path starting from $u$ to $v$, and both $u$ and $v$ are major vertices. Similarly, a directed $k$-$uv$-chain $(u(k)v]$ (or $(uv]$ for short) if $u$ is a leaf and $v$ is a major vertex.
Moreover, $(uv]$ is called {\it open} and $[uv]$ is called {\it closed}.
Note that closed $uv$-chain  has two orientations $[uv]$ and $[vu]$.

A major vertex is called {\it generalized major handle} if it is incident to at least $\Delta-1$ open chains. A chain is called a {\it terminal chain} if it is incident to a leaf. It is clear that each major handle must be a major generalized handle.

For convenience, we define $T_{vu}(v)$ as the connected component of $T-vu$ containing $v$ (normally it is call a {\it subtree of $u$}), and $T_{u}(uv)$ as the graph derived from adding $\{u,uv\}$ to $T_{vu}(v)$. Similarly,  define $T_{[vu]}(v)$ and $T_{u}([uv])$ if $[vu]$ (or $[uv]$) is a chain of $T$. A subtree $T'$ of $T$ is called a {\it strong subtree} if there is no vertex $u\in V_3(T)\cap V_2(T')$. Clearly, $T_u([uv])$ is a strong subtree of $T$ but $T_{[vu]}(v)$ is not.
%\begin{center}
%\epsfig{file=Tvuv.eps, width=9cm}\\
%{\small Subtrees $T_{[vu]}(v)$ and $T_u([uv])$.}
%\end{center}

From now on, we assume all trees $T$ are of maximum degree $3$. Let $T'$ be a tree obtained from $T$ by replacing each $uv$-chain as the edge $uv$. Thus the vertex set of $T'$ is the set of leaves and major vertices of $T$.

Let $D(T)$ be a digraph with the same  vertex  set as $T'$. If $uv\in E(T')$ is the edge corresponding to the open chain $(uv]$ in $T$, then assign a direction from $u$ to $v$.  We keep the notation $(uv]$ to denote such an arc in $D(T)$.  If $uv\in E(T')$ is the edge corresponding to the closed $k$-$uv$-chain, then duplicate it into two arcs: one is the arc from $u$ to $v$ and the other is from $v$ to $u$. Keep the notation $[u(k)v]$ (or $[uv]$) and $[v(k)u]$ (or $[vu]$), which are called {\it out-arc} and {\it in-arc} of $u$, to denote such directed arcs in $D(T)$, respectively.

Next, assign non-negative weights to each major vertex of $D(T)$ (equivalently $T$) with respect to its incident in-arc. In order to distinguish weights and labels, we use circled numbers to represent the weights. Moreover, assume $T$ contains at least two major vertices and so $D(T)$ is not isomorphic to $K_{1,3}$.

\medskip

\noindent{\bf Weight Assignment:}

\noindent Initial step: Every open oriented edge $(uv]$ gives weight \round{1} to $v$.

\noindent Main procedure: If all vertices of $D(T)$ receive three weights, then Stop. Otherwise, consider each vertex $u$ with neighborhood $\{v_1, v_2, v\}$ that has received at least two weights \round{a} and \round{b} from its in-arcs $[v_1u]$ and $[v_2u]$, respectively, and one of its out-arc $[u(k)v]$ has not assigned type. (The existence of the vertex $u$ will be proved after the following examples.) Assign the out-arc $[u(k)v]$ of $u$ a type by means of Definition~\ref{def type} stated below and give $v$ a weight according to Table~\ref{table weight}. Repeat this procedure.

\medskip

After this assignment, each major vertex of $T$ receives weights belong to $\{\round{0}, \round{1}, \round{2}, \round{3}, \round{5}, \round{6}, \round{10}, \round{15}\}$.

Define the type of closed oriented chain as follows:
\begin{definition}\label{def type}
Let $[u(k)v]$ be a closed  oriented chain and $P_1,P_2$ be the other oriented chains incident to $u$, where $k\ge 2$. Suppose $P_1$ and $P_2$ give weight \round{a} and \round{b} to $u$, respectively. If $a$ and $b$ are positive, then the chain $[u(k)v]$ is
\begin{enumerate}[(1)]
\item of type $(\round{1},k)$ if $a,b \in \{1,2,3,5\}$ and $\gcd(a,b)=1$;
\item of type $(\round{a},k)$ if $a\in \{6,10,15\}$ and $\gcd(a,b)=1$;
\item of type $(\round{d},k)$ if $\gcd(a,b)=d\in \{2,3,5\}$;
\item of type $(\round{0},k)$ if gcd$(a,b)\in\{6,10,15\}$,
\end{enumerate}
where $\gcd(a,b)$ is the greatest common divisor (gcd) of $a$ and $b$.
We also define $[u(k)v]$ of type $(\round{0},k)$ when one of $a$ and $b$ is zero.
\end{definition}

\begin{rem}\label{rem-condition4} {\rm Case~(4)  is equivalent to $a=b\in\{6,10,15\}$.
}\end{rem}

\begin{table}[H]
\begin{tabular}{|*{8}{@{\,}c}|c|}
\multicolumn{9}{c}{$[u(0)v]$ gives \round{6} to $v$;
$[u(1)v]$ gives \round{15} to $v$.}\\
\hline
\multicolumn{8}{|c|}{Types of oriented $[u(k)v]$  ($k\ge2$)} & \ Weight assigned to $v$ \ \\\hline
$(\round{1},2)$, & $(\round{1},4^+)$, & $(\round{2},7^+)$, & $(\round{3},7^+)$,  & $(\round{5},6^+)$, & $(\round{6},5^+)$, & $(\round{10},4^+)$, & $(\round{15},5^+)$ & \round{1}\\\hline
$(\round{1},3)$, & & $(\round{2},5)$, & $(\round{3},6)$, & $(\round{5},3)$, & $(\round{6},2)$, & $(\round{10},3)$, & $(\round{15},3)$ & \round{2}\\\hline
& & $(\round{2},6)$, & $(\round{3},5)$, & & & $(\round{10},2)$, & $(\round{15},4)$ & \round{3}\\\hline
& & $(\round{2},4)$, & & $(\round{5},5)$, & $(\round{6},4)$,& & $(\round{15},2)$ & \round{5}\\\hline
& & $(\round{2},2)$, & $(\round{3},3)$, & $(\round{5},4)$,  & $(\round{6},3)$ & & & \round{6}\\\hline
& & & $(\round{3},2)$ & & & & & \round{10}\\\hline
& & & $(\round{3},4)$, & $(\round{5},2)$ & & & & \round{15}\\\hline
$(\round{0},2^+)$, & & $(\round{2},3)$ &  &  & & & & \round{0}\\\hline
\multicolumn{9}{c}{$k^+$ means all the integers not less than $k$.}
\end{tabular}
\caption{The weight is given to the terminal $v$ by $[u(k)v]$, $k\ge 0$.}\label{table weight}
\end{table}

Now we explain why the algorithm can assign all notes with three weights.

Let $r=\lceil\mbox{diam}(T')/2\rceil$ be the radius of $T'$, where diam$(T')$ is the diameter of $T'$. It is known that $T'$ contains one or two centers depending on whether diam$(T')$ is even or odd. Let $c$ be a center of $T'$ and consider $T'$ as a rooted tree with root $c$.

If $r=1$, then $T'$ is $K_{1,3}$. Clearly, three weights can be assigned to the center. Hence we assume  $r\ge 2$. In this case there are at least two pairs of leaves adjacent to their fathers, respectively. After the initial step, there are at least two vertices that receive two weights from their sons. So we may perform the main procedure. By removing all leaves of $T'$, we denote the resulting tree by $T''$. There are at least two pairs of leaves of $T''$ adjacent to their fathers respectively if the radius of $T''$ is still greater than 1. For those leaves of $T''$, they have already received two weights from their sons. By the same argument as before, there is at least two vertices receiving two weights. Therefore the assignment may continue and remove the leaves repeatedly until the radius of the last tree becomes 1. In this case, the last tree is either $K_{1,3}$ or $P_2$. Note that, up to now each vertex receives weights from its two sons. Also, after removing the leaves at each iteration, the height of $T''$ decreases exactly one. As a result, the last vertex(ices) receiving weights must be the center(s).

For the first case, the main procedure implies $c$ receives three weights from its sons. For the second case, the main procedure implies $c$ and the other center $c'$ receive two weights from their sons. By the rule of the assignment, we perform the main procedure on $c$ and $c'$ and they receive weights from $[cc']$ and $[c'c]$, respectively.

Now we back to consider the original graph $D(T)$. If we remove all arcs that have types, the resulting graph is isomorphic to $T''$. At this stage, all centers receive three weights and other vertices receive two weights. Thus the main procedure can be performed from the center(s) to its(their) descendants step by step.

\begin{rem}\label{rem-subtree}{\rm From the above assignment, it is easy to see that the weight given to a vertex $u$ from its in-arc $[vu]$ is uniquely determined by the subtree of $u$ containing the vertex $v$.
}
\end{rem}

\begin{definition}\label{def bad subtree}
  A strong subtree $T^*$ of $T$ with $\Delta(T)=3$ is called {\it bad} if it satisfies the following conditions:
\begin{enumerate}[(1)]
\item $V_3(T^*)\neq \varnothing$.

\item For each generalized major handle $u$ of $T^*$, two of its
terminal chains are closed $3$-chains in $T$, or one of its incident
chains is a closed $0$-chain or $1$-chain in $T$.

\item There is no major vertex adjacent to another major vertex in $T^*$
and no 2-vertex adjacent to two different major vertices in $T^*$.

\item There is a vertex $u\in V_3(T^*)$ satisfying one of the following conditions
\begin{enumerate}[{(4.}1)]
\item One of its incident closed chain $[uv]$ is type
$(\round{2},3)$ in $T$.

\item Two of its incident chains give the same weight \round{6},
\round{10} or \round{15} to $u$.

\item Vertex $u$ receives three positive weights and the greatest common divisor of these weights is greater than $1$.
\end{enumerate}

Such vertex $u$ is called a {\it bad vertex}. A vertex is {\it good} if it is not bad.
\end{enumerate}

A tree $T$ is called a {\it bad tree} if it contains a bad strong subtree.
\end{definition}
\begin{ex} {\rm Consider the following tree $T$. The middle figure is the corresponding tree $T'$ and the right figure is the digraph $D(T)$.
%$$\begin{array}{ccc}
%\epsfig{file=gooda.eps, width=5cm}&\epsfig{file=goodb.eps, width=5cm}&\epsfig{file=good0.eps, width=5cm}\\
% T & T' & D(T)\end{array}$$
%The iterations according to the algorithm are demonstrated as follows:
%$$\begin{array}{cc}
%\epsfig{file=good1.eps, width=6cm} & \epsfig{file=good2.eps, width=6cm}\\ \epsfig{file=good3.eps, width=6cm} & \epsfig{file=good4.eps, width=6cm}
%\end{array}$$
Hence $T$ is a good tree. \rsq
}
\end{ex}
\begin{ex}{\rm Consider another tree $T$, we obtain $D(T)$ similarly.
%$$\begin{array}{ccc}
%\epsfig{file=bada.eps, width=5cm}&\epsfig{file=goodb.eps, width=5cm}&\epsfig{file=bad0.eps, width=5cm}\\
% T & T' & D(T)\end{array}$$
% After weight assignment we have
% $$\begin{array}{cc}
%\epsfig{file=good1.eps, width=6cm} & \epsfig{file=bad2.eps, width=6cm}\\ \epsfig{file=bad3.eps, width=6cm} & \epsfig{file=bad4.eps, width=6cm}
%\end{array}$$
Vertex $v$ receives the same weight $\round{6}$. The gcd of three weights of vertex $u$ received is $\round{2}$. Closed chain $[uw]$ is of type $(\round{2},3)$. Hence $u$ and $v$ are bad vertices and $T$ is a bad tree. The following figure is a bad subtree of $T$.

%\centerline{\epsfig{file=badaa.eps, width=5cm}}
In this tree, $u$ is incident with a closed 0-chain and a closed 3-chain in $T$; $v$ is incident with a closed 0-chain and an open 0-chain in $T$.}
\rsq
\end{ex}
\begin{rem}\label{rem-major}{\rm For any generalized major handle $u$, each terminal chain gives \round{1} to it. Then the closed chain, if any, is of type $(\round{1}, k)$. Thus $u$ does not satisfy (4) of Definition~\ref{def bad subtree} and it is not bad vertex. In other words, a bad vertex is incident to at least two closed chains.}\end{rem}

\begin{rem}\label{rem good weight}{\rm
In a good tree,  every oriented chain can give positive weight to its major terminal. Therefore,  \round{0} weight can only be given to its major terminal in bad tree.}\end{rem}

A configuration is called $<333>$ if it has a 3-vertex adjacent to two 3-vertices, or $<32323>$ if it has a 3-vertex adjacent two 2-vertices which is adjacent to another 3-vertex respectively (see Figure \ref{cfg 32323}). It is obvious that $<333>$ associates a bad subtree $T^*$ which is the subtree induced by the vertex set $\{u_1, u_2, u, x\}$, because both two closed 0-chain give \round{6} to their common adjacent 3-vertex. Similarly, $<32323>$ associates a bad subtree $T^*$ which is the subtree induced by the vertex set $\{u_1, u_2, y_1, y_2, u, x\}$.

%\begin{figure}[H]
%\begin{center}
%  \includegraphics[scale=0.5]{333.eps}
%  \hspace{5em}
%  \includegraphics[scale=0.5]{32323.eps}
%   \end{center}
%\caption {\label{cfg 32323} Configurations  $<333>$ and $<32323>$}
% \end{figure}

\begin{lemma}\label{lem bad cases} Let $T$ be a tree with $\Delta=3$ and $[uv]$ be a closed chain, where $u$ is a bad vertex. Let $wu$-chain $P$ and $w'u$-chain $Q$ be the other chains incident to $u$. Suppose that  $[vu]$, $P$ and $Q$ give positive weights \round{a}, \round{b} and \round{c} to $u$, respectively. Then one of the fo  statements holds:
\begin{enumerate}[$(1)$]
\item $b=c\in \{6,10,15\}$;

\item $a\in \{6,10,15\}$ and either $a\in \{b,c\}$ with $\gcd(b,c)=1$ or $\gcd(a,b,c)\in \{2,3,5\}$;

\item $a\in \{2,3,5\}$ and is a factor of $\gcd(b,c)$.
    \end{enumerate}\end{lemma}

\proof Since $[vu]$, $P$ and $Q$ give weights to their ends, they are not of type $(\round{2},3)$  according to Table~\ref{table weight}. Moreover, by Remark~\ref{rem-major}, at least one of $P$ and $Q$ is closed.

Since $u$ is a bad vertex and by Definition~\ref{def bad subtree}, we have three cases: (A) $[uv]$, $[uw]$ or $[uw']$ is of type $(\round{2},3)$, or (B) two of $a$, $b$ and $c$ are same and in $\{6,10,15\}$, or (C) gcd$(a,b,c)>1$.

\begin{enumerate}[(A)]\setlength{\parindent}{2ex}
\item By symmetry, we may assume that $[uv]$ is of type $(\round{2},3)$. Then $\gcd(b, c)=2$ by Definition~\ref{def type}. Since $[v(3)u]$ is not of type $(\round{2},3)$ and $a>0$, $a\in\{2,6\}$ by Table~\ref{table weight}. Thus, $\gcd(a,b,c)=2$. It is referred to (C).

\item If $b=c\in \{6,10,15\}$, then (1) holds. Suppose $a=b\in\{6,10,15\}$ or $a=c\in\{6,10,15\}$. Without loss of generality, we may assume $a=b\in\{6,10,15\}$ and $b\neq c$. Thus, $\gcd(a,b,c)\in\{1,2,3,5\}$. If $\gcd(a,b,c)=1$, then $\gcd(b,c)=1$. Otherwise, $\gcd(a,b,c)\in\{2,3,5\}$ and we have (2).

\item$\gcd(a,b,c)>1$. Clearly $a\ne 1$.

Suppose $a\in \{2,3,5\}$. Since $a$ is prime and $\gcd(a,b,c)>1$, we have $\gcd(a,b,c)=a$ and so $a$ is a factor of $\gcd(b,c)$, which implies (3).

Suppose $a\in \{6,10,15\}$.  If $a\in\{b,c\}$, then it is referred to (B). Let $a\notin\{b,c\}$. If $b=c\in \{6,10,15\}$, then (1) holds. Otherwise, we have $\gcd(b,c)\in \{1,2,3,5\}$. Since $\gcd(a,b,c)>1$, $\gcd(a,b,c)\in \{2,3,5\}$ and  we have (2).\qed

\end{enumerate}

\begin{lemma}
\label{lem good good}
Let $T$ be a tree with $\Delta=3$. Assume a closed chain $[u(k)v]$ and the other two chains incident to $v$ give \round{a}, \round{d} and \round{d'} to $v$ respectively; $[v(k)u]$ and the other two chains incident to $u$ give \round{b}, \round{c} and \round{c'} to $u$ respectively, where $a, b, c, c', d, d'$ are positive and $k\ge 2$. Then $u$ is bad if and only if $v$ is bad. In other words, $u$ is good if and only if $v$ is good.
\end{lemma}

\proof By symmetry, assume $u$ is bad. From Lemma~\ref{lem bad cases} we know that $b\ne 1$. Hence, by Table~\ref{table weight}, $k\le 6$ and $[v(k)u]$ is one of the type listed at the second to seventh rows of Table~\ref{table weight}.

Note that, gcd$(c,c')\in\{6,10,15\}$ implies $c=c'\in\{6,10,15\}$ and so $a=0$, which is not a case. Therefore gcd$(c,c')\in\{1,2,3,5\}$.

Suppose $k=2$.
\begin{enumerate}[({A1.}1)]\vspace*{-3mm}
\item $[v(2)u]$ is of type $(\round{2},2)$. By Definition~\ref{def type}, gcd$(d,d')=2$ and from Table~\ref{table weight}, $b=6$.

By Lemma~\ref{lem bad cases}, without loss of generality, we have $c=6$ with gcd$(c,c')=1$; or gcd$(6,c,c')=2$; or gcd$(6,c,c')=3$. From the note above and Definition~\ref{def type}, $[u(2)v]$ corresponds to type $(\round{6},2)$; $(\round{2},2)$ or $(\round{3},2)$, respectively, which gives \round{2}, \round{6} and \round{10} to $v$ by Table~\ref{table weight}. Thus gcd$(a,d,d')=2$ and so $v$ is bad.

\item $[v(2)u]$ is of type $(\round{3},2)$. By Definition~\ref{def type}, gcd$(d,d')=3$ and from Table~\ref{table weight}, $b=10$.

By Lemma~\ref{lem bad cases}, without loss of generality, we have $c=10$ with gcd$(c,c')=1$; or gcd$(10,c,c')=2$; or gcd$(10,c,c')=5$. Similar to (A1.1), $[u(2)v]$ corresponds to type $(\round{10},2)$; $(\round{2},2)$ or $(\round{5},2)$, respectively, which yields $a=3, 6,15$ by Table~\ref{table weight}. Thus gcd$(a,d,d')=3$ and so $v$ is bad.

\item  $[v(2)u]$ is of type $(\round{5},2)$. By Definition~\ref{def type}, gcd$(d,d')=5$ and from Table~\ref{table weight}, $b=15$.

By Lemma~\ref{lem bad cases}, without loss of generality, we have $c=15$ with gcd$(c,c')=1$; or gcd$(15,c,c')=3$; or gcd$(15,c,c')=5$. Similarly, $[u(2)v]$ corresponds to type $(\round{15},2)$; $(\round{3},2)$ or $(\round{5},2)$, respectively, which yields $a=5,10,15$. Then gcd$(a,d,d')=5$  and hence $v$ is bad.

\item $[v(2)u]$ is of type $(\round{6},2)$. By Definition~\ref{def type}, without loss of generality, $d=6$ with gcd$(d,d')=1$ and from Table~\ref{table weight}, $b=2$.

By Lemma~\ref{lem bad cases}, $2$ is a divisor of gcd$(c,c')$. From the note above, gcd$(c,c')=2$.
Therefore, $[u(2)v]$ is of type $(\round{2},2)$ and hence $a=6$. Since both $a$ and $d$ are 6, $v$ is bad.

\item $[v(2)u]$ is of type $(\round{10},2)$. By Definition~\ref{def type}, without loss of generality, $d=10$ with gcd$(d,d')=1$ and from Table~\ref{table weight}, $b=3$.

By Lemma~\ref{lem bad cases}, similar to the previous case, gcd$(c,c')=3$.
Therefore, $[u(2)v]$ is of type $(\round{3},2)$ and hence $a=10$. Since both $a$ and $d$ are 10, $v$ is bad.

\item $[v(2)u]$ is of type $(\round{15},2)$.
By Definition~\ref{def type}, without loss of generality, $d=15$ with gcd$(d,d')=1$ and from Table~\ref{table weight}, $b=5$.

By Lemma~\ref{lem bad cases}, similar to the previous case, gcd$(c,c')=5$.
Therefore, $[u(2)v]$ is of type $(\round{5},2)$ and hence $a=15$. Since both $a$ and $d$ are 15, $v$ is bad.
\end{enumerate}

Suppose $k=3$.
\begin{enumerate}[({A2.}1)]\vspace*{-3mm}
\item   $[v(3)u]$ is of type $(\round{3},3)$. By Definition~\ref{def type}, gcd$(d,d')=3$ and from Table~\ref{table weight}, $b=6$.

Similar to (A1.1), we have $[u(3)v]$ is of type $(\round{6},3)$, $(\round{2},3)$ or $(\round{3},3)$. Since $a>0$, $[u(3)v]$ is not of type $(\round{2},3)$. For the other cases, $a=6$ by Table~\ref{table weight}. Then gcd$(a,d,d')=3$. Hence $v$ is bad.

\item $[v(3)u]$ is of type $(\round{5},3)$ or $(\round{10},3)$. or $(\round{15},3)$ and from Table~\ref{table weight}, $b=2$.

Similar to (A1.4) we have gcd$(c,c')=2$.
Therefore, $[u(3)v]$ is of type $(\round{2},3)$ and hence it is not a case.

\item $[v(3)u]$ is of type $(\round{6},3)$. By Definition~\ref{def type}, without loss of generality, $d=6$ and gcd$(d,d')=1$ and from Table~\ref{table weight}, $b=6$.

    Similar to (A2.1) we have $a=6$, so $v$ is bad as both $a$ and $d$ are 6.
\end{enumerate}

It is similar for the cases $k=4,5,6$ and the proofs are omitted. \rsq

\medskip
\begin{lemma}
\label{lem replace}
Let $T$ and $T'$ be good trees with $\Delta=3$. Suppose that chains $[uv]$, $P$ and $Q$ give \round{a}, \round{b} and \round{c} to $v$ respectively in $T$, and $wv'$-chain gives \round{d} to $v'$ in $T'$, where $d(w)=1$ or $3$ in $T'$ and $a,b,c,d$ are positive. Let $T''$ be the tree obtained from $T$ by replacing the subtree of $v$ containing $u$ by the subtree of $v'$ containing $w$.

\begin{enumerate}[$(1)$]
\item If $d$ divides $a$, then $T''$ is also a good tree.

\item If $a\in\{\gcd(d,b),\gcd(d,c)\}$ and $a\in \{2,3,5\}$, then $T''$ is also a good tree.
\end{enumerate}
\end{lemma}%\addtocounter{theorem}{-1}

\proof Under the hypothesis and by the weighted assignment, we have the following figures:
%$$\begin{array}{ccccc}
%\epsfig{file=merge1.eps, height=3cm}&\qquad\qquad & \epsfig{file=merge2.eps, height=3cm} &\qquad\qquad & \epsfig{file=merge3.eps, height=3cm}\\
%T & & T'& & T''
%\end{array}$$
Suppose to the contrary that $T''$ is bad. Hence $T''$ has a strong subtree $T^*$ satisfies the conditions in Definition~\ref{def bad subtree}. Suppose $v$ ($v'$) is not a vertex of $T^*$.  The $T^*$ lies in one of the components of $T''-v$. Thus $T^*$ is a strong subtree of either $T$ or $T'$, which contradicts with $T$ and $T'$ being good. As a result, $v$ is a vertex of $T^*$.

\begin{enumerate}[{Step}~1:]

\item  Vertex $v$ is good in $T^*$.

Suppose vertex $v$ is bad in $T^*$.
By Lemma~\ref{lem bad cases}, there are four possible cases: (A) $b=c\in \{6,10,15\}$; (B) $d\in \{6,10,15\}$, $d\in \{b,c\}$ and $\gcd(b,c)=1$; (C) $d\in \{6,10,15\}$ and $\gcd(d,b,c)\in \{2,3,5\}$; (D) $d\in \{2,3,5\}$ and $d$ is a factor of $\gcd(b,c)$.

\begin{enumerate}[(1)]
\item $d$ divides $a$.

Since $v$ is good in $T$, gcd$(a,b,c)=1$. Moreover, if $b=c$, then $b\not\in\{6,10, 15\}$. The former property implies that gcd$(d,b,c)=1$. Combining these two properties, only (B) will occur, i.e., $d\in \{6,10,15\}$, $d\in \{b,c\}$ and $\gcd(b,c)=1$. Without loss of generality, we may assume $b=d$. Since $d$ divides $a$, we have $a=b=d$. This implies that $v$ is bad in $T$ and contradiction occurs.

\item $a\in\{\gcd(d,b),\gcd(d,c)\}$ and $a\in \{2,3,5\}$.

Without loss of generality, we assume $\gcd(d,b)=a\in \{2,3,5\}$. Since $v$ is good in $T$, gcd$(a,b,c)=1$. Since $a$ is a prime factor of $b$, gcd$(a,c)=1$. This implies that $\gcd(d,b,c)=\mbox{gcd}(a,c)=1$ and so only (B) or (D) will occur.

Suppose (B) holds. Since $d$ is composite and gcd$(d,b)=a$, which is a prime, so $d\ne b$. Hence $d=c$. This contradicts with gcd$(b,c)=1\ne a$.

Suppose (D) holds. Since $a$ and $d$ are prime, by our assumption $a=d$. Thus, gcd$(a,b,c)=1$ implies that it is impossible for $d$ being a factor of gcd$(b,c)$.

\end{enumerate}
\item Vertex $w$ is good in $T^*$.

We only need to consider when $u$ is a major vertex. Since $T'$ is a good tree, $w$ receives two positive weights from $T'$. Also, $T^*$ is bad subtree satisfying conditions of Definition~\ref{def bad subtree} implies that there are no closed 0-chain and 1-chain in $T^*$. Suppose $[vw]$ sends weight $\round{0}$ to $w$. Since we have proved that $v$ is good in $T^*$, $[vw]$ is only possible of type $(\round{0}, k)$ for some $k\ge 2$. By Remark~\ref{rem-condition4} we have $b=c\in\{6,10,15\}$, which contradicts with $v$ being good in $T^*$. Since $v$ is good, so $w$ is good in $T^*$ by Lemma~\ref{lem good good}.

\item Suppose $P$ and/or $Q$ are closed chains with other ends $x$ and $y$, respectively. Similar to Step~2, $x$ and $y$ are good in $T^*$.
\end{enumerate}

Consider $T^*$ as a rooted tree with root $v$. By the same proof as above, we can prove that the major descendants of $v$ are good. Thus $T^*$ is a good tree, which yields a contradiction. \rsq

\begin{lemma}\label{lem subtree}
If $T$ is a good tree with $\Delta=3$, then any strong subtree of $T$ is also good.
\end{lemma}

\proof Let $S$ be a strong subtree of $T$ with vertices $x_1, \dots, x_s$ that are $k$-vertices in $T$ but leaves in $S$, where $k=2,3$. Then $S$ is obtained from $T$ by suitably removing $k-1$ subtrees of each $x_i$. In order to prove this lemma, it suffices to prove that a strong subtree obtained from $T$ by removing $k-1$ subtrees of a $k$-vertex $x$ is still good, where $k=2,3$.

Let $v$ be the nearest major vertex apart from $x$ in $S$. Let $T'$ be the tree consisting of the chain $(xv]$ by adding two leaves to $v$.   Clearly $T'$ is good. Now $S$ is the tree obtained from $T$ by replacing the subtree of $v$ containing $x$ by the subtree of $v$ in $T'$ containing $x$.
By substituting $d=1$ in Lemma~\ref{lem replace}(1), we conclude that $S$ is good.

\rsq

\begin{lemma}\label{lem configuration}
 Let $T$ be a tree with $\Delta=3$ and $V_3(T)\ge 3$. If $T$ does not contain $<333>$ and $<32323>$, then $T$ contains at least one of the following configurations {\rm(}see Figure~\ref{cfg c1-4} and \ref{cfg c5-10}~{\rm)}:
\begin{enumerate}[{\rm (C1)}]
\item A path $ux_1x_2\cdots x_7$ where $u$ is a $3$-vertex and $x_1,x_2,\dots, x_7$ are $2$-vertices;

\item a leaf $u$ adjacent to a $2$-vertex $v$;

\item a $3$-vertex $v$ incident to a closed $0$-chain $[u(0)v]$ and an open $0$-chain $(y(0)v]$, where $u$ is a major handle;

\item a closed $k$-chain $[u(k)v]$, where $u$ is a major handle and
\begin{enumerate}[{\rm ({C4.}1)}]
\item  $k=2$;
\item  $k\in \{4,5,6\}$;
\end{enumerate}

\item a $3$-vertex $v$ incident to two chains: a closed $3$-chain $[u(3)v]$, where $u$ is a major handle and
\begin{enumerate}[{\rm ({C5.}1)}]
  \item an open $0$-chain $(y(0)v]$;

  \item  a closed $0$-chain $[y(0)v]$, where $y$ is a major handle;

  \item a closed $1$-chain $[y(1)v]$, where $y$ is a major handle;

  \item a closed $3$-chain $[y(3)v]$, where $y$ is a major handle;
\end{enumerate}

\item a $3$-vertex $v$ incident to three chains: a closed $1$-chain $[u(1)v]$, where $u$ is a major handle; an open $0$-chain $(y(0)v]$; and a closed $k$-chain $[v(k)w]$, where
\begin{enumerate}[{\rm ({C6.}1)}]
\item  $k=0$;
\item $k\in\{3,5,6\}$;
\end{enumerate}

\item a $3$-vertex $v$ incident to three chains: a closed $1$-chain $[u(1)v]$, where $u$ is a major handle; an open 0-chain $(y(0)v]$; and a closed $2$-chain $[v(2)w]$, while $w$ is incident to
\begin{enumerate}[{\rm ({C7.}1)}]
\item an open $0$-chain $(w'(0)w]$;

\item a closed $0$-chain $[w'(0)w]$, where $w'$ is a major handle;

\item a closed $1$-chain $[w'(1)w]$, where $w'$ is a major handle;

\item a closed $3$-chain $[w'(3)w]$, where $w'$ is a major handle;

\item a closed $2$-chain $[w'(2)w]$ and $w'$ is incident to an open $0$-chain $(w_1(0)w']$ and a closed $1$-chain $[w_2(1)w']$, where $w_2$ is a major handle;
\end{enumerate}

\item a $3$-vertex $v$ incident to three chains: a closed $1$-chain $[u(1)v]$, where $u$ is a major handle; a closed $0$-chain $[u'(0)v]$, where $u'$ is a major handle; and a closed $2$-chain $[v(2)w]$, while $w$ is incident to

\begin{enumerate}[{\rm ({C8.}1)}]
  \item an open $0$-chain $(w'(0)w]$;

  \item a closed $0$-chain $[w'(0)w]$, where $w'$ is a major handle;

  \item a closed $1$-chain $[w'(1)w]$, where $w'$ is a major handle;

  \item a closed $3$-chain $[w'(3)w]$, where $w'$ is a major handle;

  \item a closed $2$-chain $[w'(2)w]$ and $w'$ is incident to an open $0$-chain $(w_1(0)w']$ and a closed $1$-chain $[w_2(1)w']$, where $w_2$ is a major handle;

  \item a closed $2$-chain $[w'(2)w]$ and $w'$ is incident to a closed $0$-chain $[w_1(0)w']$ and a closed $1$-chain $[w_2(1)w']$, where $w_1$ and $w_2$ are major handles;
\end{enumerate}

\item a $3$-vertex $v$ incident to three chains: a closed $1$-chain $[u(1)v]$, where $u$ is a major handle; an open $0$-chain $(y(0)v]$; and a closed $4$-chain $[v(4)w]$, while $w$ is incident to

\begin{enumerate}[{\rm ({C9.}1)}]
  \item an open $0$-chain $(w'(0)w]$;

  \item a closed $0$-chain $[w'(0)w]$, where $w'$ is a major handle;

  \item a closed $1$-chain $[w'(1)w]$, where $w'$ is a major handle;

  \item a closed $3$-chain $[w'(3)w]$, where $w'$ is a major handle;

  \item a closed $2$-chain $[w'(2)w]$ and $w'$ is incident to an open $0$-chain $(w_1(0)w']$ and a closed $1$-chain $[w_2(1)w']$, where $w_2$ is a major handle;

  \item a closed $2$-chain $[w'(2)w]$ and $w'$ is incident to a closed $0$-chain $[w_1(0)w']$ and a closed $1$-chain $[w_2(1)w']$, where $w_1$ and $w_2$ are major handles;

  \item a closed $4$-chain $[w'(4)w]$ and $w'$ is incident to an open $0$-chain $(w_1(0)w']$ and a closed $1$-chain $[w_2(1)w']$, where $w_2$ is a major handle;
\end{enumerate}

\item a $3$-vertex $v$ incident to three chains: a closed $1$-chain $[u(1)v]$, where $u$ is a major handle; a closed $0$-chain $[u'(0)v]$, where $u'$ is a major handle; and a closed $k$-chain $[v(k)w]$, where $k\in\{3,4,5,6\}$.
\end{enumerate}
\end{lemma}
The proof of Lemma~\ref{lem configuration} is not difficult but tedious only. The main idea of the proof of Lemma~\ref{lem configuration} is to consider every vertex of a longest path of a tree, similar to those in \cite{WWF2006} and
\cite{ZMQ2012}. The remaining parts are only careful and tedious analysis and we omit the proof here.

\section{Main results}

\begin{theorem}[\cite{GRIGGS1992}]\label{lem two type}
For every tree $T$, $\Delta+1\le\lambda(T)\le \Delta+2$.
\end{theorem}

\begin{lemma}[\cite{WWF2006}]\label{lem 04}
If $T$ is a tree with $\Delta=3$ and $f$ is a $4$-$L(2,1)$-labeling of $T$, then $f(u)=0$ or $4$ for every major vertex $u$.
\end{lemma}

\begin{theorem}\label{th}
Let $T$ be a tree with $\Delta=3$. Then $\lambda(T)=4$ if and only if $T$ is good.
\end{theorem}

We give the proof of Theorem~\ref{th} by considering the sufficiency and necessity of $T$ is good.

\subsection{Sufficiency}
Theorem~\ref{lem two type} shows that $\lambda(T)\geq 4$ for any tree $T$ with $\Delta=3$. Hence in this subsection we assume that $T$ is a good tree. It suffices to show that $T$ has a $4$-$L(2,1)$-labeling. It is easy to obtain a $4$-$L(2,1)$-labeling if $|V_3(T)|\leq 2$, therefore we assume that $|V_3(T)|\geq 3$.

\begin{rem}\label{rem-symmetric} {\rm Suppose that a tree $T$ with $\Delta = 3$ has a 4-$L(2, 1)$-labeling $f$ using the label set $\mathcal{B}=\{0, 1, 2, 3, 4\}$. Define $f' = 4-f$, then $f'$ is also a 4-$L(2, 1)$-labeling of $T$, which is called {\it the symmetric labeling of $f$}.
}\end{rem}

We  prove that $T$ has a 4-$L(2,1)$-labeling using \bb by induction on $|T|$. Since $T$ is good, $T$ does not contain the configurations $<333>$ and $<32323>$. By Lemma~\ref{lem configuration}, we only need to deal with cases  (C1)-(C10).

\begin{enumerate}[(C1)]
\item There is a path $ux_1x_2\cdots x_7$, where $u$ is a 3-vertex and $x_1, x_2, \dots, x_7$ are 2-vertices. Assume $v$ is the another neighbor of $x_7$ besides $x_6$.

Let $T'=T\setminus \{x_2,x_3,\dots,x_6\}$. Then $T'$ consists of two components, say $T_1$ and $T_2$. Assume that $u\in V(T_1)$ and $v\in V(T_2)$. Since $V_3(T)\cap V_2(T')=\varnothing$, $T_i$ are strong subtrees of $T$, where $i=1,2$. By Lemma~\ref{lem subtree}, they are good. Thus, $T'$ has a 4-$L(2,1)$-labeling $f$ by induction hypothesis. By Remark~\ref{rem-symmetric} and Lemma~\ref{lem 04} we may assume that $f(u)=0$. In this case $f(x_1)\in\{2,3,4\}$.
 \begin{enumerate}[({1.}1)]
 \item Suppose that $f(x_1)=2$. By Remark~\ref{rem-symmetric} we may assume that $f(v)\in\{0,1,2\}$. Note that the label of $x_7$ has some restrictions depended on the label of $v$. For example, when $f(v)=1$ then $f(x_7)\in \{3,4\}$. Following is the label assignment for $x_2, x_3, x_4, x_5, x_6$:
     $$\begin{array}{|c|*{5}{|c}||c|c|}\hline
     x_1 & x_2 & x_3 & x_4 & x_5 & x_6 & x_7 & v\\\hhline{|========|}
     2 & 4 & 1 & 3 & 0 & 4 & 2 & 0\\\hline
     2 & 4 & 0 & 2 & 4 & 1 & 3 & 0\\\hline
     2 & 4 & 1 & 3 & 0 & 2 & 4 & 0\\\hline
     2 & 4 & 0 & 2 & 4 & 0 & 3 & 1\\\hline
     2 & 4 & 1 & 3 & 0 & 2 & 4 & 1\\\hline
     2 & 4 & 0 & 3 & 1 & 4 & 0 & 2\\\hline
     \end{array}$$
  \item Suppose that $f(x_1)=3$. In order to have a $4$=$L(2,1)$-labeling, the label of $x_2$ and $x_3$ must be $1$ and $4$, respectively.
       By Remark~\ref{rem-symmetric} we may assume that $f(v)\in\{0,3,2\}$. Following is the label assignment for $x_2, x_3, x_4, x_5, x_6$ except when $f(v)=0$ and $f(x_7)=3$:
     $$\begin{array}{|c|*{5}{|c}||c|c|}\hline
     x_1 & x_2 & x_3 & x_4 & x_5 & x_6 & x_7 & v\\\hhline{|========|}
     3 & 1 & 4 & 2 & 0 & 4 & 2 & 0\\\hline
     3 & 1 & 4 & 0 & 3 & 1 & 4 & 0\\\hline
     3 & 1 & 4 & 2 & 0 & 4 & 1 & 3\\\hline
     3 & 1 & 4 & 0 & 2 & 4 & 0 & 3\\\hline
     3 & 1 & 4 & 0 & 3 & 1 & 4 & 2\\\hline
     \end{array}$$
  For the case $f(x_7)=3$ and $f(v)=0$, we relabel $T_2$ by the symmetric labeling of $f$. Hence the labels of $x_7$ and $v$ are $1$ and $4$, respectively. Then we label $x_2, x_3, x_4, x_5, x_6$ by $1,4,2,0,3$, respectively.
  \item Suppose that $f(x_1)=4$. Similar to the case above by relabeling $T_2$ if necessary, we have the following assignment $x_2, x_3, x_4, x_5, x_6$:
  $$\begin{array}{|c|*{5}{|c}||c|c|}\hline
     x_1 & x_2 & x_3 & x_4 & x_5 & x_6 & x_7 & v\\\hhline{|========|}
     4 & 2 & 0 & 3 & 1 & 4 & 2 & 0\\\hline
     4 & 1 & 3 & 0 & 4 & 1 & 3 & 0\\\hline
     4 & 2 & 0 & 4 & 1 & 3 & 0 & 4\\\hline
     4 & 1 & 3 & 0 & 2 & 4 & 1 & 3\\\hline
     4 & 2 & 0 & 4 & 2 & 0 & 4 & 1\\\hline
     4 & 2 & 0 & 4 & 2 & 0 & 4 & 2\\\hline
     \end{array}$$
 \end{enumerate}
As a result, we only consider $T$ contains closed $k$-chain, where $k\le 6$, in the remaining cases.

\item There is a leaf $u$ adjacent to a 2-vertex $v$. Let $w$ be the other neighbor of $v$.
Let $T'=T-u$, then $T'$ has a 4-$L(2,1)$-labeling $f$ by
Lemma~\ref{lem subtree} and induction hypothesis. We may label $u$ by the element in $\mathcal{B} \setminus\{f(w), f(v), f(v)-1, f(v)+1\}$.

\item There is a 3-vertex $v$ incident to a closed 0-chain
$[u(0)v]$ and an open 0-chain $(y(0)v]$ such that $u$ is a major
handle. Let $y_1$ and $y_2$ be the leaves adjacent to the handle
$u$. Since $|V_3(T)|\ge3$, $v$ must be incident to another closed $k$-chain
$[v(k)w]=vx_1x_2\cdots x_kw$. From Case~(C1) and because $T$ does not contain the configuration $<333>$, so $1\le k\le 6$ in the remaining parts of this proof. Note that,
$[u(0)v]$ gives \round{6} to $v$ and $(y(0)v]$ gives \round{1} to $v$, so
$[v(k)w]$ is type $(\round{6},k)$. We consider the following cases
with different values of $k$:
\begin{enumerate}[(3.1)]
\item  $k=1$. Let $T'=T-\{y_1,y_2\}$. Then $T'$ has a 4-$L(2,1)$-labeling
$f$ by Lemma~\ref{lem subtree} and induction
hypothesis. By Lemma~\ref{lem 04} and Remark~\ref{rem-symmetric}, we may assume $f(v)=0$ and hence $f(w)=4$ and $f(x_1)=2$.  Then relabel $y$ by 3 and $u$
by 4. Finally, set $f(y_1)=1$ and $f(y_2)=2$ in $T$.
\item $k=2$. Then $[v(2)w]$ is of type $(\round{6},2)$ and gives \round{2} to $w$.
We construct a new tree $T'=T_{wx_2}(w)+[w(3)z_1]+[z_1(0)z_2)+[z_1(0)z_3)$ and $z_0\in
[w(3)z_1]$ be the neighbor of $w$. Note that it is isomorphic to $T-y$. It is easy to see that
$[z_1(3)w]$ also gives \round{2} to $w$. Then $T'$ is a good tree and has a
4-$L(2,1)$-labeling $f$ by Lemma~\ref{lem replace} and
induction hypothesis. We may assume $f(w)=0$ and $f(z_0)=\alpha$. Hence $\alpha\in\{3,4\}$, otherwise $z_1$ cannot be labeled under $f$ by Lemma~\ref{lem 04}. Now we label $x_2$ by $\alpha$, i.e., either 3 or 4.
Hence, assign proper label sequence $03140$ or $04204$ to $wx_2x_1vu$ in
$T$. Thus, $f$ can be extended to $T$ after labeling
the leaves adjacent to $u$ and $v$ easily since
$f(u),f(v)\in\{0,4\}$.

\item $k=3$. Then $[v(3)w]$ is of type $(\round{6},3)$ and gives \round{6} to $w$.
Let $T'=T_{wx_3}(w)+[w(0)z_1]+[z_1(0)z_2)+[z_1(0)z_3)$. It is easy to  see that $[z_1(0)w]$ gives \round{6} to $w$. Therefore, $T'$ is a good tree by Lemma~\ref{lem replace} and has a 4-$L(2,1)$-labeling $f$ by induction hypothesis. By Lemma~\ref{lem 04}, we may assume $f(w)=0$ and $f(z_1)=4$. Hence, assign proper label sequence $041304$ to $wx_3x_2x_1vu$ in $T$. Thus, $f$ can be extended to $T$ after labeling the leaves adjacent to $u$ and $v$.

\item $k=4$. Then $[v(4)w]$ is of type $(\round{6},4)$ and gives \round{5} to $w$.
Let $T'=T_{wx_4}(w)+[w(2)z_1]+[z_1(0)z_2)+[z_1(1)z_3]+[z_3(0)z_4)+[z_3(0)z_5)$ and $z_0\in [w(2)z_1]$ be the neighbor of $w$. It is easy to see that $[z_1(2)w]$ is of type $(\round{15},2)$ and gives \round{5} to $w$. Therefore, $T'$ is a good tree by Lemma~\ref{lem replace} and has a 4-$L(2,1)$-labeling $f$ by induction hypothesis.
It is clear that $f(w),f(z_1),f(z_3)\in\{0,4\}$ by Lemma~\ref{lem 04}. If $f(z_0)\in\{0,4\}$, then the vertex adjacent with $z_1$ and $z_3$ and the vertex adjacent with $z_1$ and $w$ must be labeled by $2$ which is impossible. So we get $f(z_0)\notin \{0,4\}$. Assume $f(w)=0$. We label $x_4$ by $f(z_0)\in\{2,3\}$. The possible label sequence of the path $wx_4x_3x_2x_1vu$ is $0241304$ or $0314204$ in $T$. Thus, $f$ can be extended to $T$ after labeling the leaves adjacent to $u$ and $v$.

\item $k=5$.  Let $T'=T_{x_5x_4}(x_5)$. So $T'\subset T$ is a good tree by Lemma~\ref{lem subtree}. By induction hypothesis, $T'$ has a 4-$L(2,1)$-labeling $f$ with $f(w)=0$. Then $f(x_5)=2$, 3, or 4. According to these cases, we label $x_5x_4x_3x_2x_1vu$ in $T$ with sequence
$2403140$, $3140240$, or $4130240$, respectively.
\item $k=6$. Similar to the above case, $T'=T_{x_5x_4}(x_5)$ has a 4-$L(2,1)$-labeling $f$ with $f(w)=4$ by induction hypothesis.

If $f(x_6)=0$, then $f(x_5)=2$, 3 or 4 which is the same as (3.5).

If $f(x_6)=1$, then and $f(x_5)=3$ or $4$. We label the path $x_5x_4x_3x_2x_1vu$ in $T$ with label sequence $3041304$ or $4204204$ accordingly.

If $f(x_6)=2$, then $f(x_5)=0$. We label the path $x_5x_4x_3x_2x_1vu$ in $T$ with label sequence $0314204$.

Thus, $f$ can be extended to $T$ after labeling the leaves adjacent to $u$ and $v$.
\end{enumerate}

\item There is a 3-vertex $v$ incident to a closed $k$-chain $[u(k)v]=ux_1x_2\cdots x_kv$ such that $u$ is a major handle, where $k\in\{2,4,5,6\}$.
\begin{enumerate}[(4.1)]
\item $k=2$.
 Let $T'=T_{x_2x_1}(x_2)$. Then $T'\subset T$ is a good tree by Lemma~\ref{lem subtree}. By the induction hypothesis, $T'$ has a 4-$L(2,1)$-labeling $f$ with $f(v)=0$. Whatever the label of $x_2$ is in $T'$, we always can assign proper label sequence $0240$, $0314$ or $0420$ to the path $vx_2x_1u$ in $T$. Thus, $f$ can be extended to $T$ after labeling the leaves adjacent to $u$.

\item $k\in \{4,5,6\}$.
Let $T'=T_{x_4x_3}(x_4)$. Then $T'\subset T$ is a good tree by Lemma~\ref{lem subtree}. By the induction hypothesis, $T'$ has a 4-$L(2,1)$-labeling $f$. By Remark~\ref{rem-symmetric} we may assume that $f(x_5)\in\{0,1,2\}$ when $k=5,6$ (or $f(v)=0$ when $k=4$).

If $f(x_5)=0$ (or $f(v)=0$ when $k=4$), then $f(x_4)=2$, 3  or 4. We label the path $x_4x_3x_2x_1u$ in $T$ with label sequence $24130$, $31420$ or $41304$ accordingly.

If $f(x_5)=1$, then $f(x_4)=3$ or 4. We label the path $x_4x_3x_2x_1u$ in $T$ with label sequence $30420$ or $40314$ accordingly.

If $f(x_5)=2$, then $f(x_4)=0$ or 4.  We label the path $x_4x_3x_2x_1u$ in $T$ with label sequence $03140$ or $41304$ accordingly.

Thus, $f$ can be extended to $T$ after labeling the leaves adjacent to $u$.
\end{enumerate}
\item
There is a 3-vertex $v$ incident to two chains: one is a closed 3-chain $[u(3)v]=ux_1x_2x_3v$ such that $u$ is a major handle. The other chain $Q$ is an open $0$-chain $(y(0)v]$ or a closed chain $[y(k)v]$ with $k=0,1,3$ such that $y$ is a major handle.
\begin{enumerate}[(5.1)]
\item  $Q=(y(0)v]$.
Let $T'=T_{x_3x_2}(x_3)$. By the induction hypothesis, $T'$ has a 4-L(2,1)-labeling $f$ with $f(v)=0$. If $f(x_3)=2$, exchange the labels of $y$ and $x_3$. Hence, $f(x_3)=3$ or $4$. We can assign $03140$ or $04204$ to $vx_3x_2x_1u$ in $T$. Thus, $f$ can be extended to $T$ after labeling the leaves adjacent to $u$.

\item $Q=[y(0)v]$ such that $y$ is a major handle.
Let $P$ be the chain incident to $v$ besides $[u(3)v]$ and $[y(0)v]$. If $P$ is open, then we can label all vertices of $T$ easily. Assume $P=uy_1\cdots y_kw$ is closed. Since $T$ does not contain $<333>$ and by Case~(C1), $1\le k\le 6$. Note that, $[u(3)v]$ and $[y(0)v]$ give \round{2} and \round{6} to $v$ respectively, so $[v(k)w]$ is of type $(\round{2},k)$. Since $T$ is good, $k\ne 3$. Next we consider the following cases depending on the values of $k$:

\begin{enumerate}[(5.2-1)]
\item
 $k=1$.
Let $T'=T_{x_3x_2}(x_3)$. Then $T'$ has a 4-$L(2,1)$-labeling $f$ with $f(v)=0$ by Lemma~\ref{lem subtree} and the induction hypothesis. Then $f(y_1)=2$ and $f(y)=4$ and it deduces $f(x_3)=3$ in $T$. Assign proper label sequence $03140$ to $vx_3x_2x_1u$. Thus, $f$ can be extended to $T$ after labeling the leaves adjacent to $u$ and $y$.

\item $k=2$. Then $[v(2)w]$ is of type $(\round{2},2)$ and gives \round{6} to $w$.
Let $T'=T_{wy_2}(w)+[w(0)z_1]+[z_1(0)z_2)+[z_1(0)z_3)$. Same as (3.3) we have $[z_1(0)w]$ gives \round{6} to $w$, $T'$ is a good tree and has a 4-$L(2,1)$-labeling $f$ with $f(w)=0$ and $f(z_1)=4$. We assign proper label sequence $04203140$ to $wy_2y_1vx_3x_2x_1u$ and 4 to $y$.  Thus, $f$ can be extended to $T$ after labeling the leaves adjacent to $u$ and $y$.

\item $k=4$. Then $[v(4)w]$ is of type $(\round{2},4)$ and gives \round{5} to $w$.
Let $T'=T_{wy_4}(w)+[w(2)z_1]+[z_1(0)z_2)+[z_1(1)z_3]+[z_3(0)z_4)+[z_3(0)z_5)$ and $z_0\in [w(2)z_1]$ be the neighbor of $w$. Same as (3.4) we have $[z_1(2)w]$ gives \round{5} to $w$, $T'$ is a good tree and has a 4-$L(2,1)$-labeling $f$ with $f(w)=0$ and $f(z_0)\in \{2,3\}$. We assign proper label sequence $0240241304$ or $0314203140$ to $wy_4y_3y_2y_1vx_3x_2x_1u$ in $T$ and label $y$ with 0 or 4, respectively. Thus, $f$ can be extended to $T$ after labeling the leaves adjacent to $u$ and $y$.

\item $k=5$. Then $[v(5)w]$ is of type $(\round{2},5)$ and gives \round{2} to $w$.
Let $T'=T_{wy_5}(w)+[w(3)z_1]+[z_1(0)z_2)+[z_1(0)z_3)$ and $z_0\in [w(3)z_1]$ be the neighbor of $w$. Same as (3.2) we have $[z_1(3)w]$ gives \round{2} to $w$, $T'$ is a good tree and has  a 4-$L(2,1)$-labeling $f$ with $f(w)=0$ and $f(z_0)=3$ or 4. Hence, we assign proper label sequence $0314024$ or $0413024$ to $wy_5y_4y_3y_2y_1v$ and then let $f(y)=0$ and assign $41304$ to $vx_3x_2x_1u$ in $T$. Thus $f$ can be extended to $T$ after labeling the leaves adjacent to $u$ and $y$.

\item $k=6$. Then $[v(6)w]$ is of type $(\round{2},6)$ and gives \round{3} to $w$.
Let $T'= T_{wy_6}(w)+[w(4)z_1]+[z_1(0)z_2)+[z_1(1)z_3]+[z_3(0)z_4)+[z_3(0)z_5)$ and $z_0\in [w(4)z_1]$ be the neighbor of $w$. It is easy to see that $[z_1(4)w]$ is of type by $(\round{15},4)$ and gives \round{3} to $w$. So $T'$ is a good tree by Lemma~\ref{lem replace} and has a 4-$L(2,1)$-labeling $f$ with $f(w)=0$ by induction hypothesis. Since $f(z_1), f(z_3)\in\{0,4\}$, we can check that $f(z_0)\in\{2,4\}$. We assign proper label sequence $024130241304$ or $041304203140$ to $wy_6\cdots y_1vx_3x_2x_1u$  and label $y$ with 0 or 4 in $T$, respectively. Thus, $f$ can be extended to $T$ after labeling the leaves adjacent to $u$ and $y$.
\end{enumerate}

\item $Q=[y(1)v]$ such that $y$ is a major handle.
Let $T'=T_{x_3x_2}(x_3)$. Then $T'$ has a 4-$L(2,1)$-labeling $f$ with $f(v)=0$ by the induction hypothesis. Hence $f(y)=4$. Let $z$ be the common neighbor of $v$ and $y$. Then $f(z)=2$ and so  $f(x_3)=3$ or 4. We assign proper sequence $03140$ or $04204$ to $vx_3x_2x_1u$, accordingly. Thus, $f$ can be extended to $T$ after labeling the leaves adjacent to $u$.

\item $Q=[y(3)v]$ such that $y$ is a major handle. Let $Q=yy_1y_2y_3v$.
Note that both $[u(3)v]$ and $[y(3)v]$ give \round{2} to $v$. Let $T'=T_{vy_3}(v)+[v(0)z_1]+[z_1(0)z_2)+[z_1(0)z_3)$. It is easy to see that $[z_1(0)v]$ gives \round{6} to $v$. By the Lemma~\ref{lem replace} and the induction hypothesis, $T'$ is a good tree and it has a 4-$L(2,1)$-labeling $f$ with $f(v)=0$. Hence $f(z_1)=4$ and so we properly label $y_3$ with $4$ in $T$. Next, we assign proper label sequence $04204$ to $vy_3y_2y_1y$. Finally, $f$ can be extended to $T$ after labeling the leaves adjacent to $y$.

\end{enumerate}

\item There is a 3-vertex $v$ incident to three chains: a closed 1-chain $[u(1)v]=uxv$ such that $u$ is a major handle; an open 0-chain $(y(0)v]$; and a closed $k$-chain $[v(k)w]=vy_1y_2\cdots y_kw$, where $k\in\{0,3,5,6\}$.

\begin{enumerate}[(6.1)]

\item $k=0$.

Let $T'=T_{xu}(x)$. Then $T'\subset T$ is a good tree by Lemma~\ref{lem subtree}. By induction hypothesis, $T'$ has a 4-$L(2,1)$-labeling $f$ with $f(v)=0$. Hence $f(w)=4$. Relabel $y$ by 3 and $x$ by 2. Finally, label $u$ by 4 in $T$. Thus, $f$ can be extended to $T$ after labeling the leaves adjacent to $u$.

\item $k=3$. Let $T'=T_{xu}(x)$. Then $T'\subset T$ is a good tree by Lemma~\ref{lem subtree} and has 4-$L(2,1)$-labeling $f$ with $f(v)=0$ by induction hypothesis. It is easy to see that $f(y_1)\neq 2$ since $f(v),f(w)\in \{0,4\}$. Then label $x$ with 2. Finally, label $u$ with 4 in $T$. Thus, $f$ can be extended to $T$ after labeling the leaves adjacent to $u$ and $v$.

\item $k=5$ or $6$. Let $T'=T_{y_ky_{k-1}}(y_k)$. Then $T'\subset T$ is a good tree by Lemma~\ref{lem subtree} and has 4-$L(2,1)$-labeling $f$ with $f(w)=0$ by induction hypothesis.  No matter $f(x_k)=2$, $3$ or $4$, we always assign proper label sequence to $wy_k\cdots y_1v$ in $T$ by $0240240$, $0314204$ or $0420314$ when $k=5$; assign $02403140$, $03140314$ or $04130240$ when $k=6$. For each case, $f(v)\in \{0,4\}$ and $f(y_1)\neq 2$. As a result, we can set $f(x)=2$ and $f(u)=4-f(v)$. Thus, $f$ can be extended to $T$ after labeling the leaves adjacent to $u$ and $v$.
\end{enumerate}
\item There is a 3-vertex $v$ incident to three chains: a closed 1-chain $[u(1)v]=uu_1v$ where $u$ is a major handle; an open 0-chain $(y(0)v]$; and a closed 2-chain $[v(2)w]=vx_1x_2w$.
\begin{enumerate}[(7.1)]
\item $w$ is incident to an open 0-chain $(w'(0)w]$.
Let $T'=T_{x_2x_1}(x_2)$. Then $T'\subset T$ is a good tree by Lemma~\ref{lem subtree} and has 4-$L(2,1)$-labeling $f$ with $f(w)=0$. If $f(x_2)=4$, then exchange the labels of $w'$ and $x_2$. As a result, $f(x_2)\in \{2,3\}$. Assign proper label sequence $024024$ or $031402$ to $wx_2x_1vu_1u$ accordingly. Thus, $f$ can be easily extended to $T$.

\item $w$ is incident to a closed 0-chain $[w'(0)w]$ such that $w'$ is a major handle.
Let $T'=T_{x_2x_1}(x_2)$. Then $T'$ is a good tree and hence has 4-$L(2,1)$-labeling $f$ with $f(w)=0$ by induction hypothesis. Thus, $f(w')=4$. Similar to the (7.1), $f$ can be extended to $T$.

\item $w$ is incident to a closed 1-chain $[w'(1)w]$ such that $w'$ is a major handle.

Let $[w(k)w'']=wy_1y_2\cdots y_kw''$ be the chain incident to $w$ besides $[v(2)w]$ and $[w'(1)w]$. Since $T$ does not contain $<32323>$, $k\ne 1$ and thus $k\in\{0,2,3,4,5,6\}$. Note that, $[v(2)w]$ is of type $(\round{15},2)$ and gives \round{5} to $w$; while $[w'(1)w]$ gives \round{15} to $w$. Hence $[w(k)w'']$ is of type $(\round{5},k)$. We consider the following cases with different values of $k$:
\begin{enumerate}[(7.3-1)]
\item $k=0$.
Let $T'=T_{x_2x_1}(x_2)$. Then $T'$ is a good tree and hence has 4-$L(2,1)$-labeling $f$ with $f(w)=0$ by induction hypothesis. Then $f(w'')=4$ and $f(w')=4$, which induces $f(x_2)=3$. Assign proper label sequence $1420$ to $x_1vu_1u$. Thus, $f$ can be extended to $T$.

\item $k=2$. Note that $[w(2)w'']$ is of type $(\round{5},2)$ and gives \round{15} to $w''$.
Let $T'=T_{w''y_2}(w'')+[w''(1)z_1]+[z_1(0)z_2)+[z_1(0)z_3)$, where $[w''(1)z_1]=w''z_0z_1$. It is easy to see that $[z_1(1)w'']$ gives \round{15} to $w''$ too. Therefore, $T'$ is a good tree and has 4-$L(2,1)$-labeling $f$ with $f(w'')=0$ by Lemma~\ref{lem replace} and induction hypothesis. Then $f(z_1)=4$ and $f(z_0)=2$. Hence, assign proper label sequence $0240$ to $w''y_2y_1w$ in $T$. Note that $f(y_1)=4$ and so, similar to (7.3-1), $f$ can be extended to $T$.

\item $k=3$. Note that $[w(3)w'']$ is of type $(\round{5},3)$ and gives \round{2} to $w''$. Let $T'=T_{w''y_3}(w'')+[w''(3)z_1]+[z_1(0)z_2)+[z_1(0)z_3)$ and $z_0\in [w''(3)z_1]$ be the neighbor of $w''$. Same as (3.2) we have $[z_1(3)w'']$ gives \round{2} to $w''$, $T'$ is a good tree and has a 4-$L(2,1)$-labeling $f$ with $f(w'')=0$ and $f(z_0)=3$ or 4.
So we may label $y_3$ by $f(z_0)$. And assign proper label sequence $3140$ or $4204$ to $y_3y_2y_1w$ in $T$ according to $f(z_0)=3$ or $4$. Note that $f(y_1)\in \{0,4\}$ and so, similar to (7.3-1), $f$ can be extended to $T$.

\item $k=4$. Note that $[w(4)w'']$ is of type $(\round{5},4)$ and gives \round{6} to $w''$.
Let $T'=T_{w''y_4}(w'')+[w''(0)z_1]+[z_1(0)z_2)+[z_1(0)z_3)$. Same as (3.3) we have $[z_1(0)w'']$ gives \round{6} to $w''$, $T'$ is a good tree and has a 4-$L(2,1)$-labeling $f$ with $f(w'')=0$ and $f(z_1)=4$. Then we can assign proper sequence $41304$ to $y_4y_3y_2y_1w$ in $T$. Note that $f(y_1)=0$ and so, similar to (7.3-3), $f$ can be extended to $T$.
\item $k=5$. Note that $[w(5)w'']$ is of type $(\round{5},4)$ and gives \round{5} to $w''$.
Let $T'=T_{y_5w''}(w'')+[w''(2)z_1]+[z_1(0)z_2)+[z_1(1)z_3]+[z_3(0)z_4)+[z_3(0)z_5)$ and $z_0\in [w''(2)z_1]$ be the neighbor of $w''$. Same as (3.4) we have $[z_1(2)w'']$ gives \round{5} to $w''$, $T'$ is a good tree and has a 4-$L(2,1)$-labeling $f$ with $f(w'')=0$ and $f(z_0)\in \{2,3\}$. Then we assign proper label sequence $241304$ or $314204$ to $y_5\cdots y_1w$ in $T$ according to $f(z_0)=2$ or 3. Note that $f(y_1)=0$ and so, similar to (7.3-3), $f$ can be extended to $T$.

\item $k=6$.
Let $T'=T_{y_6y_5}(y_6)$. Then $T'$ is a good tree and hence has a 4-$L(2,1)$-labeling $f$ with $f(w'')=0$ by induction hypothesis. Whatever the label of $y_6$ is, assign proper label sequence $2403140$, $3140240$ or $4130240$ to $y_6\cdots y_1w$ in $T$. Therefore, $f(y_1)=4$ and hence $f$, similar to (7.3-1), $f$ can be extended to $T$.

\end{enumerate}
\item $w$ is incident to a closed 3-chain $[w'(3)w]=w'y_1y_2y_3w$, where $w'$ is a major handle.
Let $T'=T_{y_3y_2}(y_3)$. Then $T'$ is a good tree and hence has a 4-$L(2,1)$-labeling $f$ with $f(w)=0$ by induction hypothesis.
Since $f(u),f(v)\in \{0,4\}$ and $f(w)=0$, we may show that $f(x_2)=2$ or $3$. If $f(y_3)=2$, then $f(x_2)=3$. Exchange the labels of $x_2$ and $y_3$. Next, relabel $x_1vu_1u$ by the proper label sequence $4024$, relabel $y$ by 3, and relabel the leaves adjacent to $u$ by 0 and 1. As a result, $f(y_3)\neq 2$. Next, assign proper label sequence $3140$ or $4204$ to $y_3y_2y_1w'$ according to $f(y_3)=3$ or 4. Finally, $f$ can be extended to $T$ after labeling the leaves adjacent to $w'$.

\item $w$ is incident to a closed 2-chain $[w'(2)w]=w'y_1y_2w$ such that $w'$ is incident to an open 0-chain $(w_1(0)w']$ and a closed 1-chain $[w_2(1)w']$, where $w_2$ is a major handle. Note that both $[v(2)w]$ and $[w'(2)w]$ is of type $(\round{15},2)$ and give \round{5} to $w$, respectively. Let $T'=T_{wx_2}(w)+[w(1)z_1]+[z_1(0)z_2)+[z_1(0)z_3)$, where $[w(1)z_1]=wz_0z_1$.
Then $[z_1(1)w]$ gives \round{15} to $v$. By Lemma~\ref{lem replace} and the induction hypothesis, $T'$ is a good tree and it has a 4-$L(2,1)$-labeling $f$ with $f(w)=0$. Hence $f(z_1)=4$ and $f(z_0)=2$. Now we may label $x_2$ by 2 and assign proper label sequence $24024$ to $x_2x_1vu_1u$ in $T$. Thus, $f$ can be extended to $T$.
\end{enumerate}
\item There is a 3-vertex $v$ incident to three chains: a closed 1-chain $[u(1)v]=uu_1v$ such that $u$ is a major handle; a closed 0-chain $[u'(0)v]$ such that $u'$ is a major handle; a closed $2$-chain $[v(2)w]=vx_1x_2w$.
\begin{enumerate}[(8.1)]

\item $w$ is incident to an open 0-chain $(w'(0)w]$. Let $P$ be the chain incident to $v$ besides $[v(2)w]$ and $(w'(0)w]$. If $P$ is open, then we can label all vertices of $T$ easily. Assume $P=[w(k)w'']=wy_1\cdots y_kw''$ is closed.
Note that, $[u(1)v]$ and $[u'(0)v]$ give \round{15} and \round{6} to $v$ respectively. Hence,  $[v(2)w]$ is of type $(\round{3},2)$ and it gives \round{10} to $w$. Furthermore, $[w(k)w'']$ is of type $(\round{10},k)$. We consider the following cases with different values of $k$:
\begin{enumerate}[({8.1}-1)]

\item $k=0$.
Let $T'=T_{x_2x_1}(x_2)$. Then $T'$ is a good tree and hence has a 4-$L(2,1)$-labeling $f$ with $f(w)=0$ by induction hypothesis. Therefore, we have $f(w'')=4$. Relabel $x_2$ with 3 and $w'$ with 2. Assign proper label sequence $31420$ to $x_2x_1vu_1u$ in $T$ and let $f(u')=0$. Thus, $f$ can be extended to $T$.

\item $k=1$.
Let $T'=T_{x_2x_1}(x_2)$. Then $T'$ is a good tree and hence has 4-$L(2,1)$-labeling $f$ with $f(w)=0$ by induction hypothesis. Thus, $f(w'')=4$ and $f(y_1)=2$. Relabel $x_2$ with 3 and $w'$ with 4 and so, similar to (8.1-1), $f$ can be extended to $T$.

\item $k=2$. Note that $[w(2)w'']$ is of type $(\round{10},2)$ and gives \round{3} to $w$.
Let $T'=T_{w''y_2}(w'')+[w''(4)z_1]+[z_1(0)z_2)+[z_1(1)z_3]+[z_3(0)z_4)+[z_3(0)z_5)$ and $z_0\in [w''(4)z_1]$ be the neighbor of $w''$. Similar to (5.2-5), we have $f(z_0)\in\{2,4\}$ and hence we may let $f(y_2)=f(z_0)$. Then we assign proper label sequence $240$ or $420$ to $y_2y_1w$ in $T$. Similar to (8.1-1), we can label $x_2$ with 3. As a result, $f$ can be extended to $T$.

\item $k=3$. Then $[w(3)w'']$ is of type $(\round{10},3)$ and gives \round{2} to $w$.

Let $T'=T_{y_3w''}(w'')+[w''(3)z_1]+[z_1(0)z_2)+[z_1(0)z_3)$ and $z_0\in [w''(3)z_1]$ be the neighbor of $w''$.
Same as (3.2) we have $T'$ is a good tree and has a 4-$L(2,1)$-labeling $f$ with $f(w'')=0$ and $f(z_0)=3$ or 4.
We assign proper label sequence $3140$ or $4204$ to $y_3y_2y_1w$ in $T$. Consequently, we can label $u'$ with 2 and $x_2$ with 1 or 3 and so, similar to (8.1-1), $f$ can be extended to $T$.

\item $k=4$, 5 or 6.

Let $T'=T_{y_4y_3}(y_4)$. Then $T'$ is a good tree and hence has a 4-$L(2,1)$-labeling $f$ with $f(y_5)\in\{0,1,2\}$ when $k=5,6$ and $f(w'')=0$ when $k=4$ by induction hypothesis.

If $f(y_5)=0$ (or $f(w'')=0$ when $k=4$), then $f(y_4)=2$, 3 or 4. We label the path $y_4y_3y_2y_1w$ in $T$ with label sequence $24024$, $31420$ or $41304$ accordingly.

If $f(y_5)=1$, then $f(y_4)=3$ or 4. We label the path $y_4y_3y_2y_1w$ in $T$ with label sequence $30420$ or $40240$ accordingly.

If $f(y_5)=2$, then $f(y_4)=0$ or 4. We label the path $y_4y_3y_2y_1w$ in $T$ with label sequence $03140$ or $41304$ accordingly.

For each case, $f(w)\in \{0,4\}$ and $f(y_1)\notin\{1,3\}$. Therefore, similar to (8.1-1), $f$ can be extended to $T$.
\end{enumerate}
\item $w$ is incident to a closed 0-chain $[w'(0)w]$, where $w'$ is a major handle. Recall that $[v(2)w]$ gives \round{10} to $w$ in $T$.
Let $T'=T_{wx_2}(w)+[w(3)z_1]+[z_1(0)z_2)+[z_1(0)z_3)$ and $z_0\in [w(3)z_1]$ be the neighbor of $w$. Clearly $[z_1(3)w]$ gives \round{2} to $w$. By  Lemma~\ref{lem replace} and the induction hypothesis, $T'$ is a good tree and has a 4-$L(2,1)$-labeling $f$ with $f(w)=0$ and hence $f(w')=4$. Moreover, $f(z_0)\in \{1,3\}$ and we label $x_2$ with $f(z_0)$ in $T$. Thus, similar to (8.1-1), $f$ can be extended to $T$.

\item  $w$ is incident to a closed 1-chain $[w'(1)w]=w'y_1w$, where $w'$ is a major handle.
Let $T'=T_{wx_2}(w)+[w(2)z_1]+[z_1(0)z_2)+[z_1(1)z_3]+[z_3(0)z_4)+[z_3(0)z_5)$ and $z_0\in [w(3)z_1]$ be the neighbor of $w$. Recall  that $[v(2)w]$ gives \round{10} to $w$ in $T$. Clearly,  $[z_1(2)w]$ gives \round{5} to $w$. By  Lemma~\ref{lem replace} and the induction hypothesis, $T'$ is a good tree and has a 4-$L(2,1)$-labeling $f$ with $f(w)=0$ and hence $f(w')=4$,  $f(y_1)=2$ and $f(z_0)\in \{1,3\}$. Similar to (8.2), $f$ can be extended to $T$.
\item $w$ is incident to a closed 3-chain $[w'(3)w]=w'y_1y_2y_3w$, where $w'$ is a major handle.
Let $T'=T_{wx_2}(w)+[w(3)z_1]+[z_1(0)z_2)+[z_1(0)z_3)$ and $z_0\in [w(3)z_1]$ be the neighbor of $w$. Note that it is  same as the tree $T'$ in (C8.2). Hence, $T'$ is a good tree and has a 4-$L(2,1)$-labeling $f$ with $f(w)=0$.
Since $f(w),f(w'),f(z_1)\in \{0,4\}$, $f(z_0), f(y_3)\in \{3,4\}$. By Remark~\ref{rem-symmetric} we may exchange the labels of $y_3$ and $z_0$ if $f(y_3)=3$. As a result, $f(y_3)=4$. Reassign proper sequence $04204$ to $wy_3y_2y_1w'$ and relabel the leaves of $w'$. Then $x_2$ may be assigned with 3 in $T$ and, similar to (8.1-1), $f$ can be extended to $T$.
\item $w$ is incident to a closed 2-chain $[w'(2)w]=w'y_1y_2w$ and $w'$ is incident to an open 0-chain $(w_1(0)w']$ and a closed 1-chain $[w_2(1)w']=w_2w''w'$, where $w_2$ is a major handle. Let
$T'=T_{wx_2}(w)+[w(2)z_1]+[z_1(0)z_2)+[z_1(1)z_3]+[z_3(0)z_4)+[z_3(0)z_5)$ and $z_0\in [w(2)z_1]$ be the neighbor of $w$. Since $[z_1(3)w]$ gives \round{5} to $w$, by  Lemma~\ref{lem replace} and the induction hypothesis, $T'$ is a good tree and has a 4-$L(2,1)$-labeling $f$ with $f(w)=0$.
Since $f(w)=0$ and $f(w'),f(w_2),f(z_1),f(z_3)\in \{0,4\}$,  $f(z_0),f(y_2)\in \{2,3\}$. By symmetry we may exchange the labels of $y_2$ and $z_0$ such that $f(y_2)=2$ and $f(z_0)=3$ if necessary. Assign proper sequence $024024$ to $wy_2y_1w'w''w_2$, and relabel the leaves of $w'$ and $w_2$. Label $x_2$ with 3 in $T$ and so $f$ can be extended to $T$ same as (8.1-1).

\item $w$ is incident to a closed 2-chain $[w'(2)w]$ and $w'$ is incident to a closed 0-chain $[w_1(0)w']$ and a closed 1-chain $[w_2(1)w']$, where $w_1$ and $w_2$ are major handles.
Note that $[v(2)w]$ gives \round{5} to $w$. On the other hand $[w_2(1)w']$ and $[w_1(0)w']$ give \round{15} and \round{6} to $v$, respectively. Then $[w'(2)w]$ is of type $(\round{3},2)$ and gives \round{10} to $w$. Therefore, $w$ is a bad vertex which is not a case.

\end{enumerate}

\item There is a 3-vertex $v$ incident to three chains: a closed 1-chain $[u(1)v]$, where $u$ is a major handle; an open 0-chain $(y(0)v]$; and a closed 4-chain $[v(4)w]=vx_1x_2x_3x_4w$. Note that $[v(4)w]$ is of type $(\round{15},4)$ and gives \round{3} to $w$.
\begin{enumerate}[(9.1)]
\item $w$ is incident to an open 0-chain $(w'(0)w]$.
Let $T'=T_{x_4x_3}(x_4)$. Then $T'$ is a good tree and hence has a 4-$L(2,1)$-labeling $f$ with $f(w)=0$ by induction hypothesis. Exchange the labels of $x_4$ and $w'$ if $f(x_4)=3$. As a result, $f(x_4)\in \{2,4\}$. Assign proper label sequence $2413024$ or $4130420$ to $x_4x_3x_2x_1vu_1u$. Thus, $f$ can be extended to $T$.

\item $w$ is incident to a closed 0-chain $[w'(0)w]$, where $w'$ is a major handle.
Recall that $[v(4)w]$ gives \round{3} to $w$.
Let $T'=T_{wx_4}(w)+[w(1)z_1]+[z_1(0)z_2)+[z_1(0)z_3)$ and $z_0\in [w(1)z_1]$ be the neighbor of $w$. Note that $[z_1(1)w]$ gives \round{15} to $w$ in $T'$. Hence $T'$ is a good tree by Lemma~\ref{lem replace} and has a 4-$L(2,1)$-labeling $f$ with $f(w)=0$ by induction hypothesis.  Therefore, $f(z_1)=4$ and $f(z_0)=2$ and we label $x_4$ by $f(z_0)=2$. Same as (9.1), $f$ can be extended to $T$.

\item $w$ is incident to a closed 1-chain $[w'(1)w]$, where $w'$ is a major handle.
Let $T'=T_{wx_4}(w)+[w(0)z_1]+[z_1(0)z_2)+[z_1(0)z_3)$. Note that $[z_1(0)w]$ gives \round{6} to $w$ in $T'$. So $T'$ is a good tree by Lemma~\ref{lem replace} and has a 4-$L(2,1)$-labeling $f$ with $f(w)=0$ by induction hypothesis.
Then $f(z_1)=4$. Hence, we may label $x_4$ by $f(z_1)=4$. As the same as (9.1), $f$ can be extended to $T$.

\item $w$ is incident to a closed 3-chain $[w'(3)w]=w'y_1y_2y_3w$, where $w'$ is a major handle.
Let $T'=T_{x_4x_3}(x_4)$. Then $T'$ is a good tree and hence has a 4-$L(2,1)$-labeling $f$ with $f(w)=0$ by induction hypothesis. Whatever the labels of $x_4$ and $y_3$ are, it is easy to find a relabeling strategy such that $f(x_4)\in \{2,4\}$ and $f(y_3)\in \{3,4\}$. Same as (9.1), $f$ can be extended to $T$.

\item $w$ is incident to a closed 2-chain $[w'(2)w]=w'y_1y_2w$ and $w'$ is incident to an open 0-chain $(w_1(0)w']$ and a closed 1-chain $[w_2(1)w']$, where $w_2$ is a major handle. Let $T'=T_{x_4x_3}(x_4)$. Then $T'$ is a good tree and has a 4-$L(2,1)$-labeling $f$ with $f(w)=0$ by induction hypothesis. Whatever the labels of $x_4$ and $y_2$ are, it is easy to find a relabeling strategy such that $f(x_4)\in \{2,4\}$ and $f(y_2)\in \{2,3\}$. Same as (9.1), $f$ can be extended to $T$.
\item $w$ is incident to a closed 2-chain $[w'(2)w]=w'y_1y_2w$ and $w'$ is incident to a closed 0-chain $[w_1(0)w']$ and a closed 1-chain $[w_2(1)w']=w_2zw'$, where $w_1$ and $w_2$ are major handles.
Let $T'=T_{x_4x_3}(x_4)$. Then $T'$ is a good tree and has a 4-$L(2,1)$-labeling $f$ with $f(w)=0$ by induction hypothesis. Since $\{f(w'),f(w_1)\}=\{0,4\}$ and $f(z)=2$, $f(y_1)\in \{1,3\}$. Hence $f(y_2)\ne 2$. Moreover, $f(w')\in \{0,4\}$ and $f(y_1)\ne 2$, $f(y_2)\ne 4$. As a result, $f(y_2)=3$ and so $f(x_4)\in \{2,4\}$. Same as (9.1), $f$ can be extended to $T$.

\item $w$ is incident to a closed 4-chain $[w'(4)w]$ and $w'$ is incident to an open 0-chain $(w_1(0)w']$ and a closed 1-chain $[w_2(1)w']$, where $w_2$ is a major handle. Let $T'=T_{x_4w}(w)+[w(0)z_1]+[z_1(0)z_2)+[z_1(0)z_3)$. Note that $[z_1(0)w]$ gives \round{6} to $w$ in $T'$. So $T'$ is a good tree by Lemma~\ref{lem replace} and has a 4-$L(2,1)$-labeling $f$ with $f(w)=0$ by induction hypothesis. Thus, $f(z_1)=4$ and we  label $x_4$ by $f(z_1)=4$.  Same as (9.1), $f$ can be extended to $T$.

\end{enumerate}
\item There is a 3-vertex $v$ incident to three chains: a closed 1-chain $[u(1)v]=uu_1v$, where $u$ is a major handle; a closed 0-chain $[u'(0)v]$, where $u'$ is a major handle; and a closed $k$-chain $[v(k)w]=vx_1x_2\cdots x_kw$, where $k\in\{3,4,5,6\}$.
    Note that, $[u(1)v]$ and $[u'(0)v]$ give \round{15} and \round{6} to $v$, respectively. So $[v(k)w]$ is of type $(\round{3},k)$. We consider the following cases with different values of $k$:
\begin{enumerate}[{(10}.1)]
\item  $k=3$. Then $[v(3)w]$ gives \round{6} to $w$.
Let $T'=T_{wx_3}(w)+[w(0)z_1]+[z_1(0)z_2)+[z_1(0)z_3)$. Same as (3.3) we have $f(w)=0$ and $f(z_1)=4$, where $f$ is a 4-$L(2,1)$-labeling of $T'$. We assign proper label sequence $413024$ to $x_3x_2x_1vu_1u$ and let $f(u')=4$ in $T$. Thus, $f$ can be extended to $T$.

\item $k=4$. Then $[v(4)w]$ gives \round{15} to $w$.
Let $T'=T_{wx_4}(w)+[w(1)z_1]+[z_1(0)z_2)+[z_1(0)z_3)$ and $z_0\in[w(1)z_1]$ be the neighbor of $w$. Same as (7.3-2), $T'$ has a 4-$L(2,1)$-labeling $f$  with $f(w)=0$ and hence $f(z_0)=2$. We assign proper label sequence $2413024$ to $x_4x_3x_2x_1vu_1u$ and let $f(u')=4$ in $T$. Thus, $f$ can be extended to $T$.

\item $k=5$. Then $[v(5)w]$ gives \round{3} to $w$.
Let $T'=T_{wx_5}(w)+[w(4)z_1]+[z_1(0)z_2)+[z_1(1)z_3]+[z_3(0)z_4)+[z_3(0)z_5)$ Same as (5.2-5) we obtain that $T'$ has a  4-$L(2,1)$-labeling $f$ with $f(w)=0$ and $f(z_0)=2$ or $4$. Hence, assign proper label sequence $24031420$ or $42031420$ to $x_5\cdots x_1vu_1u$ and let $f(u')=4$ in $T$. Thus, $f$ can be extended to $T$.

\item $k=6$. Then $[v(6)w]$ gives \round{2} to $w$.
Let $T'=T_{x_6w}(w)+[w(3)z_1]+[z_1(0)z_2)+[z_1(0)z_3)$ and $z_0\in [w(3)z_1]$ be the neighbor of $w$. Same as (3.2), $T'$ has  a 4-$L(2,1)$-labeling $f$ with $f(w)=0$ and $f(z_0)=3$ or 4.
Hence, assign proper label sequence $314031420$ or $420413024$ to $x_6\cdots x_1v u_1u$ and let $f(u')=0$ or $4$ in $T$ accordingly. Thus, $f$ can be extended to $T$.
\qed
\end{enumerate}
\end{enumerate}

\subsection{Necessity}

Before proving the necessity, we explain the meaning of the weights given to the vertices.

Let $K$ be a subtree of a tree $T$ with $\Delta=3$ containing a major vertex $u$. For any $w \in V(K)$, let $$S^u_{K}(w)=\{f(w)\;|\;\mbox{$f$ is a $4$-$L(2,1)$-labeling of $K$ such that $f(u)=0$}\}.$$
Clearly $S^u_{K}(w)\subseteq \{2,3,4\}$ if $w$ is a neighbor of $u$.
If we required $f(u)=4$, then we define
$$\overline{S}^u_{K}(w)=\{f(w)\;|\;\mbox{$f$ is a $4$-$L(2,1)$-labeling of $K$ such that $f(u)=4$}\}.$$ By symmetry, $\overline{S}^u_{K}(w)=\{4-a\;|\; a\in S^u_{K}(w)\}$.

\medskip

Let $T$ be a tree with $\Delta=3$ and $v\in V_3(T)$. For any $uv$-chain, let $w_1$ be the neighbor of $v$ in the $uv$-chain. Suppose $T$ has a $4$-$L(2,1)$ labeling. By symmetry we assume that the label of $v$ is $0$. Let $K=T_v(vw_1)$. We define the following rules:
{\it \begin{enumerate}[$(1)$]

\item $uv$-chain gives \round{1} to $v$ means that $S^v_{K}(w_1)=\{2,3,4\}$;

\item $uv$-chain gives \round{2} to $v$ means that $S^v_{K}(w_1)=\{3,4\}$;

\item $uv$-chain gives  \round{3} to $v$ means that $S^v_{K}(w_1)=\{2,4\}$;

\item $uv$-chain gives  \round{5} to $v$ means that $S^v_{K}(w_1)=\{2,3\}$;

\item $uv$-chain gives \round{6} to $v$ means that $S^v_{K}(w_1)=\{4\}$;

\item $uv$-chain gives  \round{10} to $v$ means that $S^v_{K}(w_1)=\{3\}$;

\item $uv$-chain gives  \round{15} to $v$ means that $S^v_{K}(w_1)=\{2\}$.
\end{enumerate}
}

\begin{rem} {\rm For Rule~(1), since $S^v_{K}(w_1)\subseteq \{2,3,4\}$, it suffices to show that for each $a\in \{2,3,4\}$, there is a $4$-$L(2,1)$-labeling $f$ such that $f(w_1)=a$ and $f(v)=0$.
For Rules~(5), (6) and (7), since we assume that $T$  has a $4$-$L(2,1)$-labeling, $S^v_{K}(w_1)\ne\varnothing$. Therefore, it suffices to show that $g(w_1)=4$, 3, and 2, respectively for any $4$-$L(2,1)$-labeling $g$ of $T$ when $g(v)=0$.
}
\end{rem}

From the definition, a weight \round{a} given to a major vertex $v$ from the closed chain $[uv]$ depends on the length of $[uv]$ and the weights \round{b} and \round{c} given to $u$ from the other two chains incident to $u$. It seems a `binary operation' ($\round{b}*\round{c}=\round{a}$) on a special system. Now we define the meaning of such weights. Hence we need to show that this meaning is still closed (well-defined) according to the `binary operation'.

\proof

If $u$ is a leaf in $T_v(vw_1)$, then $(uv]$ gives \round{1} to $v$ by Table~\ref{table weight}. For each $a\in \{2,3,4\}$, it is easy to see that there exists a $4$-$L(2,1)$ labeling $f$ for $T_v(vw_1)$ such that $f(w_1)=a$. This satisfies Rule~(1).

Now we assume $u$ is a major vertex and let $[u(k)v]=uw_kw_{k-1}\cdots w_1v$ for $k\ge 0$.

When $k=0$. Then $w_1=u$ and $[u(0)v]$ gives \round{6} to $v$. For any $4$-$L(2,1)$ labeling $g$, $g(v)=0$ implies $g(w_1)=4$. Thus it satisfies Rule~(5).

When $k=1$. $[u(1)v]$ gives \round{15} to $v$. For any $4$-$L(2,1)$ labeling $g$ with $g(v)=0$, we deduce $g(w_1)=2$ and $g(u)=4$ that satisfies Rule~(7).

Now we assume $k\ge 2$. Let $P_1$ and $P_2$ be the chains incident to $u$ besides $[uv]$, $x_1\in P_1$ and $x_2\in P_2$ be the neighbors of $u$ besides $w_k$, and $P_1$ and $P_2$ giving weight \round{y_1} and \round{y_2} to $u$, respectively.

\medskip
\noindent{\bf Properties on labels of $w_k$:}

Now we consider the subtree $H=T_{w_k}(w_ku)$.
\begin{enumerate}[(P1)]
\item If $[u(k)v]$ is of type $(\round{1},k)$, then $y_1,y_2 \in \{1,2,3,5\}$ and gcd$(y_1,y_2)=1$.
We want to show that $S^u_H(w_k)=\{2,3,4\}$.

From Rules~(1)-(4) we have $S^u_{T_1}(x_1)\cup S^u_{T_2}(x_2)=\{2,3,4\}$, where $T_1=T_u(ux_1)$ and $T_2=T_u(ux_2)$. Note that, if $y_1=y_2=1$, then by Rule~(1), $S^u_{T_1}(x_1)=S^u_{T_2}(x_2)=\{2,3,4\}$. Otherwise $S^u_{T_1}(x_1)\ne S^u_{T_2}(x_2)$.

Let $a\in\{2,3,4\}$.
Without loss of generality, we may assume $a\in S^u_{T_1}(x_1)$. It is easy to see that $|S^u_{T_1}(x_1)|\ge 2$ and $|S^u_{T_2}(x_2)|\ge 2$. So we may choose $b\in S^u_{T_1}(x_1)\setminus\{a\}$. Since either $S^u_{T_1}(x_1)=S^u_{T_2}(x_2)=\{2,3,4\}$ or $S^u_{T_1}(x_1)\ne S^u_{T_2}(x_2)$, there exists $c\in S^u_{T_2}(x_2)\setminus\{a,b\}$. According to our rules, there are $4$-$L(2,1)$-labelings $g_1$ of $T_1$ and $g_2$ of $T_2$ such that $g_1(x_1)=b$ and $g_2(x_2)=c$. Note that $g_1(u)=g_2(u)=0$. Label $w_k$ by $a$ in $H$ together with $g_1$ and $g_2$ deduce that $S^u_H(w_k)=\{2,3,4\}$.

\item If $[u(k)v]$ is of type $(\round{2},k)$, then $y_1$ and $y_2$ share the unique common prime divisor $2$. By Rules~(2), (5) and (6), we have $g_1(x_1),g_2(x_2)\neq 2$  for any $4$-$L(2,1)$-labelings $g_1$ of $T_u(ux_1)$ and $g_2$ of $T_u(ux_2)$ with $g_1(u)=g_2(u)=0$. Thus, $w_k$ can only be labeled by $2$ in $H$ and so $S^u_H(w_k)=\{2\}$. Similarly, we can show that $S^u_H(w_k)=\{3\}$ and $S^u_H(w_k)=\{4\}$ if $[u(k)v]$ is of type $(\round{3},k)$ and $(\round{5},k)$, respectively.

\item If $[u(k)v]$ is of type $(\round{6},k)$, without loss of generality,  assume $y_1=6$ and $y_2\in \{1,5\}$. By Rule~(5), $g_1(x_1)=4$ for any $4$-$L(2,1)$-labeling $g_1$ of $T_u(ux_1)$. Hence $S^u_H(w_k)\subseteq \{2,3\}$.

    For any $a\in \{2,3\}$. By Rules~(1) and (4), $\{2,3\}\subseteq S^u_{T_u(ux_2)}(x_2)$. There is a $4$-$L(2,1)$-labeling $g_2$ of $T_u(ux_2)$ such that $g_2(x_2)=b$ and $g_2(u)=0$, where  $b\in S^u_{T_u(ux_2)}(x_2)\setminus\{a\}$. Now we label $w_k$ by $a$ in $H$ and so $\{2,3\}\subseteq S^u_H(w_k)$. Therefore, $S^u_H(w_k)=\{2,3\}$.
    Similarly, we can show that $S^u_H(w_k)=\{2,4\}$ and $S^u_H(w_k)=\{3,4\}$ if $[u(k)v]$ is of type $(\round{10},k)$ and $(\round{15},k)$, respectively.
\end{enumerate}

Finally we only need to verify every type of chain in Table \ref{table weight} by using the above properties.
\begin{enumerate}[(1)]
\item $v$ receives \round{1}. We only need to show that for each $a\in\{2,3,4\}$, there is a $4$-$L(2,1)$-labeling for $T_v(vw_1)$ such that $f(w_1)=a$ and $f(v)=0$.

$\bullet$ $[u(k)v]$ is of type $(\round{1},2)$. According to the prescribed labels of $w_1$, we label $vw_1w_2u$ by $0240$, $0314$ and $0420$, respectively. By (P1), since $S^u_H(w_2)=\{2,3,4\}$, where $H=T_{w_2}(w_2u)$, we have a required labeling for the subtree $T_v(vw_1)$.

$\bullet$ $[u(k)v]$ is of type $(\round{1},4^+)$. We first assume that $k\ge 5$. No matter what is the label of $w_1$ has been assigned, there are at most twelve possible cases for the labels of $w_{k-4}w_{k-3}$. Namely, $02$, $03$, $04$, $13$, $14$, $20$, $24$, $30$, $31$, $40$, $41$ and $42$. The following six labelings and their symmetric labelings for $w_{k-4}w_{k-3}w_{k-2}w_{k-1}w_ku$ cover all those case: $024130$, $031420$, $041304$, $130420$, $140314$, $204130$ and the result follows. When $k=4$, only the first three cases are needed to consider and we obtain the result similar to the case above.

$\bullet$ Other types described at the first row of Table~\ref{table weight}, which are not mentioned here, can be verified by similar method using (P2) or (P3).
\item $v$ receives \round{2}. We need to show that $S^v_{T_v(vw_1)}(w_1)=\{3,4\}$.

$\bullet$ $[u(k)v]$ is of type $(\round{1},3)$. We label $vw_1w_2w_3u$ by $03140$ and $04130$ accordingly.  By (P1), we have a required labeling.
On the other hand, if there is a $4$-$L(2,1)$-labeling $g$ such that $g(w_1)=2$ and $g(v)=0$, then $g(w_2)=4$ and $g(u)=0$. As a result, $g(w_3)$ cannot be defined.

$\bullet$ $[u(k)v]$ is of type $(\round{10},3)$. We label $vw_1w_2w_3u$ by $03140$ and $04204$ accordingly.  By (P3), we have a required labeling. On the other hand, if there is a $4$-$L(2,1)$-labeling $g$ such that $g(w_1)=2$ and $g(v)=0$, then it is same as the case above.

$\bullet$ $[u(k)v]$ is of type $(\round{2},5)$. We label $vw_1w_2w_3w_4w_5u$ by $0314024$ and $0420420$. By (P2), we have a required labeling. On the other hand, if there is a $4$-$L(2,1)$-labeling $g$ such that $g(w_1)=2$ and $g(v)=0$, then $g(w_2)=4$, $(g(w_3),g(w_4)) =(0,2)$, $(0,3)$ or $(1,3)$, which contradicts (P2).

$\bullet$ Other types described at the second row of Table~\ref{table weight}, which are not mentioned here, can be verified by similar method using (P2) or (P3).

\item $v$ receives \round{3}. We need to show that $S^v_{T_v(vw_1)}(w_1)=\{2,4\}$.

$\bullet$ $[u(k)v]$ is of type $(\round{2},6)$. We label $vw_1w_2w_3w_4w_5w_6u$ by $02413024$ and $04130420$ respectively. By (P2), we have a required labeling. On the other hand, if there is a $4$-$L(2,1)$-labeling $g$ such that $g(w_1)=3$ and $g(v)=0$, then $g(w_2)=1$, $g(w_3)=4$ and $g(w_4)\in\{0,2\}$. But (P2) implies $g(w_6)=2$ and so $g(w_5)$ cannot be defined.

$\bullet$ $[u(k)v]$ is of type $(\round{10},2)$. We label $vw_1w_2u$ by $0240$ and $0420$ respectively.  By (P3), we have a required labeling. On the other hand, if there is a $4$-$L(2,1)$-labeling $g$ such that $g(w_1)=3$ and $g(v)=0$, then $g(w_2)=1$ and $g(u)=4$. However, (P3) implies $g(w_2)\in\{2,0\}$ and contradiction occurs.

$\bullet$ $[u(k)v]$ is of type $(\round{15},4)$. We label $vw_1w_2w_3w_4u$ by $024130$ and $041304$ respectively. By (P3), we have a required labeling. On the other hand, if there is a $4$-$L(2,1)$-labeling $g$ such that $g(w_1)=3$ and $g(v)=0$, then $g(w_2)=1$, $g(w_3)=4$ and $g(w_4)\in\{0,2\}$. Since $g(w_3)=4$, $g(u)=0$ and hence $g(w_4)=2$, which contradicts (P3).

$\bullet$ $[u(k)v]$ is of type $(\round{3},5)$. We label $vw_1w_2w_3w_4w_5u$ by $0240314$ and $0420314$ respectively.  By (P2), we have a required labeling. On the other hand, if there is a $4$-$L(2,1)$-labeling $g$ such that $g(w_1)=3$ and $g(v)=0$, then $g(w_2)=1$, $g(w_3)=4$ and $g(w_4)\in\{0,2\}$. By (P2), $(g(w_5), g(u))=(3,0)$ or $(1,4)$. This is a contradiction.

\item $v$ receives \round{5}. We need to show that $S^v_{T_v(vw_1)}(w_1)=\{2, 3\}$.

$\bullet$ $[u(k)v]$ is of type $(\round{5},5)$. We label $vw_1w_2w_3w_4w_5u$ by $0241304$ and $0314204$ respectively. By (P2), we have a required labeling. On the other hand, if there is a $4$-$L(2,1)$-labeling $g$ such that $g(w_1)=4$ and $g(v)=0$, then $(g(w_2), g(w_3), g(w_4))=(1,3,0)$, $(2,0,3)$ or $(2,0,4)$. According to (P2), $(g(w_5), g(u))=(4,0)$ or $(0,4)$ and so $g$ does not exist.

$\bullet$ $[u(k)v]$ is of type $(\round{2},4)$. We label $vw_1w_2w_3w_4u$ by $024024$ and $031420$ respectively. By (P2), we have a required labeling. On the other hand, if there is a $4$-$L(2,1)$-labeling $g$ such that $g(w_1)=4$ and $g(v)=0$, then $(g(w_2), g(w_3), g(w_4))=(1,3,0)$, $(2,0,3)$ or $(2,0,4)$. As a result, $g(u)$ cannot be defined.

$\bullet$ $[u(k)v]$ is of type $(\round{6},4)$. We label $vw_1w_2w_3w_4u$ by $024130$ and $031420$ respectively. By (P3), we have a required labeling. On the other hand, if there is a $4$-$L(2,1)$-labeling $g$ such that $g(w_1)=4$ and $g(v)=0$, then it is the same the case above.

$\bullet$ $[u(k)v]$ is of type $(\round{15},2)$.  We label $vw_1w_2u$ by $0240$ and $0314$ respectively. By (P3), we have a required labeling. On the other hand, if there is a $4$-$L(2,1)$-labeling $g$ such that $g(w_1)=4$ and $g(v)=0$, then $g(u)=4$ and $g(w_2)=2$, which contradicts (P3).

\item $v$ receives \round{6}. We only need to show that $g(w_1)=4$ for any $4$-$L(2,1)$-labeling $g$ of $T$ with $g(v)=0$.

$\bullet$ $[u(k)v]$ is of type $(\round{2},2)$. By (P2), $g(w_2)=2$ and hence $g(w_1)=4$.

$\bullet$ $[u(k)v]$ is of type $(\round{3},3)$. By (P2), the label sequence of $uw_3w_2$ under $g$ is $031$ or $413$. Therefore, $g(w_1)=4$.

$\bullet$ $[u(k)v]$ is of type $(\round{5},4)$. By (P2), the label sequence of $uw_4w_3w_2$ is $0413$, $0420$, $4031$ or $4024$. It is easy to check that only the third case exists when $g(v)=0$. It implies that $g(w_1)=4$.

$\bullet$ $[u(k)v]$ is of type $(\round{6},3)$. By (P3), the label sequence of $uw_3w_2$ is $031$, $024$, $420$ or $413$. It is easy to check that only the first case exists when $g(v)=0$. Hence $g(w_1)=4$.

\item $v$ receives \round{10}. We only need to show that $g(w_1)=3$, for any $4$-$L(2,1)$-labeling $g$ of $T$ with $g(v)=0$.

$\bullet$ $[u(k)v]$ is of type $(\round{3},2)$. By (P2), $(g(u), g(w_2))=(0,3)$ or $(4,1)$.  Since $g(v)=0$, only the last case exists and hence $g(w_1)=3$.

\item $v$ receives \round{15}. We only need to show that $g(w_1)=2$, for any $4$-$L(2,1)$-labeling $g$ of $T$ with $g(v)=0$.

$\bullet$ $[u(k)v]$ is of type $(\round{3},4)$. By (P2), the label sequence of $uw_4w_3w_2$ is $0314$ or $4130$. Since $g(v)=0$, only the first case exists and hence $g(w_1)=2$.

$\bullet$ $[u(k)v]$ is of type $(\round{5},2)$.
By (P2), the label sequence of $uw_2$ is $04$ or $40$. Since $g(v)=0$, only the first case exists and hence $g(w_1)=2$. \qed
\end{enumerate}

\medskip
\noindent {\bf The proof of necessity.}

Let $f$ be a 4-$L(2,1)$-labeling  of $T$  using label set $\mathcal{B}=\{0,1,2,3,4\}$.  Assume the contrary that $T$ is a bad tree, that is, $T$ contains a bad subtree $T^*$ such that $T^*$ satisfies the conditions in Definition \ref{def bad subtree}. Let $u$ be the bad vertex of $T^*$. Note that $u$ is a major vertex, hence we assume that $f(u)=0$. Consider the following cases depending on the reason of $u$ being the bad vertex.
\begin{enumerate}[(1)]
\item There is a closed chain $[u(3)v]=ux_1x_2x_3v$ such that $[u(3)v]$ is of type $(\round{2},3)$.

Let $P_1$ and $P_2$ be the chains incident to $u$ besides $[uv]$ and let $w_1\in P_1$ and $w_2\in P_2$ be the neighbors of $u$. Since $[u(3)v]$ is of type $(\round{2},3)$, the greatest common divisor of the weights from $P_1$ and $P_2$ to $u$ is 2. Thus, $f(w_1), f(w_2)\neq 2$ by Rules~(2), (5) and (6). This implies that $f(x_1)=2$ and $f(x_2)=4$ and hence $f(v)=0$. Thus $f$ is  not a $4$-$L(2,1)$-labeling of $T$ because $f(x_3)$ cannot be defined, contradiction occurs.

\item There are two chains incident to $u$ giving the same weight \round{6}, \round{10} or \round{15} to $u$.

Let $P_1$ and $P_2$ be these chains and $w_1\in P_1$ and $w_2\in P_2$ be the neighbors of $u$.

If $P_1$ and $P_2$ give the same weight \round{6} to $u$, then $f(w_1),f(w_2)\in \{0,4\}$ by Rule~(5). However, $f(u)=0$ which yields a contradiction.

If $P_1$ and $P_2$ give the same weight \round{10} to $u$, then $f(w_1),f(w_2)\in \{1,3\}$ by Rule~(6). However, $f(u)=0$ which yields a contradiction.

If $P_1$ and $P_2$ give the same weight \round{15} to $u$, then $f(w_1)=f(w_2)=2$ by Rule~(7). This is a contradiction.

\item The weights from all three chains incident to $u$ have the greatest common divisor 2, 3 or 5.

Let $P_1$, $P_2$ and $P_3$ be these chains and $w_1\in P_1$, $w_2\in P_2$ and $w_3\in P_3$ be the neighbors of $u$. Let $d$ be the greatest common divisor of these three weights that $u$ receives.

According to Rules~(2) to (7), if $d=2$, then $f(w_1),f(w_2),f(w_3)\in \{0,1,3,4\}$; if $d=3$, then $f(w_1),f(w_2),f(w_3)\in \{0,2,4\}$; if $d=5$, then $f(w_1),f(w_2),f(w_3)\in \{1,2,3\}$. But $f(u)=0$, which yields a contradiction.
 \end{enumerate}
This completes the proof.  \qed

\newpage
\noindent{\bf Appendix}
%\begin{figure}[H]
%
%\begin{center}
%  \includegraphics[scale=0.5]{c1-4.eps}
%  \hspace{5em}
%  \caption{\label{cfg c1-4} Configurations C1-C4}
%   \end{center}
% \end{figure}
%
% \begin{figure}[H]
%
%\begin{center}
%  \includegraphics[scale=0.5]{c5-10.eps}
%  \hspace{5em}
%  \caption{\label{cfg c5-10} Configurations C5-C10}
%   \end{center}
% \end{figure}


\begin{thebibliography}{99}

\bibitem{CHANG1996} G.J. Chang, D. Kuo, The $L(2, 1)$-labelling problem on graphs, SIAM J. Discrete Math. 9 (1996) 309-316.

\bibitem{CHANG2000}G.J. Chang, W.-T. Ke, D. Kuo, D.D.-F. Liu, R.K.Yeh, On $L(d, 1)$-labelling of graphs, Discrete Math. 220 (2000) 57-66.
\bibitem{GEORGES2000} J.P. Georges, D.W. Mauro, M.I. Stein, Labeling products of complete graphs with a condition at distance two, SIAM J. Discrete Math. 14 (2000) 28-35.

\bibitem{GEORGES1994} J.P. Georges, D.W. Mauro, M.A. Whittlesey, Relating path coverings to vertex labellings with a condition at distance two, Discrete Math. 135 (1994) 103-111.

\bibitem{GONCALVES2005} D. Gon\c{c}alves, On the $L(p,1)$-labelling of graphs, DMTCS Proceedings Volume AE (2005) 81-86.

\bibitem{GRIGGS1992} J.R. Griggs, R.K.Yeh, Labelling graphs with a condition at distance 2, SIAM J. Discrete Math. 5 (1992) 586-595.

\bibitem{HALE1980}W.K. Hale,  Frequency assignment: Theory and application, Proc. IEEE 68 (1980) 1497-1514.

\bibitem{HEUVEL2003}  J. van den Heuvel, S. McGuinness, Coloring the square of a planar graph, J. Graph Theory 42 (2003) 110-124.

\bibitem{KRAL2003} D. Kr\'{a}l, R. \v{S}krekovski, A theorem about the channel assignment problem, SIAM J. Discrete Math. 16 (2003) 426-437.

\bibitem{SAKAI1994}  D. Sakai, Labelling chordal graphs: distance two condition, SIAM J. Discrete Math. 7 (1994) 133-140.


\bibitem{WWF2006}W.F. Wang, The $L(2, 1)$-labelling of trees, Discrete Appl. Math. 154 (2006) 598-603.

\bibitem{WHITTLESEY1995}  M.A. Whittlesey, J.P. Georges, D.W. Mauro, On the $\lambda$-number of $Q_n$ and related graphs, SIAM J. Discrete Math. 8 (1995) 499-506.

\bibitem{ZMQ2012}M.Q. Zhai, C.H. Lu, J.L. Shu, A note on $L(2,1)$-labelling of Trees, Acta Math. Appl. Sin. 28 (2012) 395-400.

\end{thebibliography}
\end{document}